\documentclass[leqno,10pt]{article}
\usepackage{amsmath,amstext, amsthm, amssymb}

\usepackage{hyperref}

\usepackage[totalheight=22 true cm, totalwidth=14 true cm]{geometry}

\usepackage[all]{xy}

\newtheorem{thmintro}{Theorem}
\newtheorem{corintro}[thmintro]{Corollary}
\newtheorem{thm}{Theorem}[section]
\newtheorem{cor}[thm]{Corollary}
\newtheorem{lemma}[thm]{Lemma}
\newtheorem{prop}[thm]{Proposition}
\newtheorem{defn}[thm]{Definition}
\theoremstyle{remark}
\theoremstyle{definition}
\newtheorem{parag}[thm]{}
\newtheorem{rmk}[thm]{Remark}

\newtheorem{assumption}[thm]{Assumption}

\def\beq{\begin{equation}}
\def\eeq{\end{equation}}

\def\crash#1{}

\def\N{{\mathbb N}}

\def\C{{\mathbb C}}
\def\H{{\mathbb H}}
\def\E{{\mathbb E}}

\def\DS{{\mathbb D\!\mathbb S}}

\def\l{\left}
\def\r{\right}
\def\[[{\l[\l[}
\def\]]{\r]\r]}
\def\p{\prime}

\def\sgq{\sigma_q}

\def\dq{d_q}

\def\ord{{\rm ord}}

\def\cf{\emph{cf. }}
\def\ie{\emph{i.e. }}
\def\ds{\displaystyle}

\def\cB{{\mathcal B}}

\def\cE{{\mathcal E}}

\def\cN{{\mathcal N}}

\def\cL{{\mathcal L}}

\def\cS{{\mathcal S}}

\def\wtilde{\widetilde}

\def\veps{\varepsilon}
\def\a{\alpha}
\def\be{\beta}

\def\la{\lambda}

\author{
Lucia Di Vizio\footnote{The first author is partially supported by
the ANR contract ANR-06-JCJC-0028}\\
\small Institut de Math\'{e}matiques de Jussieu,
\small Topologie et g\'{e}om\'{e}trie alg\'{e}briques, Case 7012,\\
\small 2, place Jussieu, 75251 Paris Cedex 05, France.\\
\small e-mail: {\tt divizio@math.jussieu.fr}\\
\and\\
Changgui Zhang\\
\small Laboratoire P. PAINLEVE, U.F.R. de Math\'{e}matiques Pures et Appliqu\'{e}es, USTL\\
\small Cit\'e scientifique, 59655 Villeneuve d'Ascq CEDEX, France.\\
\small e-mail: {\tt Changgui.Zhang@math.univ-lille1.fr}}

\date{~}
\title{On $q$-summation and confluence}

\begin{document}
%\makeatletter\c@page=197\makeatother

\maketitle

\begin{abstract}
This paper is divided in two parts.
In the first part we consider a convergent $q$-analog of the
divergent Euler series, with $q\in(0,1)$, and we show how the Borel
sum of a generic Gevrey formal solution to a differential
equation can be uniformly approximated on a
convenient sector by a meromorphic solution of a corresponding $q$-difference equation.
In the second part, we work under the
assumption $q\in (1,+\infty)$. In this case, at least four
different $q$-Borel sums of a divergent power series solution of
an irregular singular analytic $q$-difference equations are spread
in the literature: under convenient assumptions we clarify the
relations among them.
\end{abstract}

\tableofcontents

%%%%%%%%%%%%%%%%%%%%%%%%%%%%%%%%%%%%%%%%%%%%%%%%%%%%%%%%%%%%%%%%%%%%
%%%%%%%%%%%%%%%%%%%%%%%%%%%%%%%%%%%%%%%%%%%%%%%%%%%%%%%%%%%%%%%%%%%%
%%%%%%%%%%%%%%%%%%%%%%%%%%%%%%%%%%%%%%%%%%%%%%%%%%%%%%%%%%%%%%%%%%%%
%%%%%%%%%%%%%%%%%%%%%%%%%%%%%%%%%%%%%%%%%%%%%%%%%%%%%%%%%%%%%%%%%%%%
%%%%%%%%%%%%%%%%%%%%%%%%%%%%%%%%%%%%%%%%%%%%%%%%%%%%%%%%%%%%%%%%%%%%
%%%%%%%%%%%%%%%%%%%%%%%%%%%%%%%%%%%%%%%%%%%%%%%%%%%%%%%%%%%%%%%%%%%%

\section*{Introduction}
\addcontentsline{toc}{section}{Introduction}

%%%%%%%%%%%%%%%%%%%%%%%%%%%%%%%%%%%%%%%%%%%%%%%%%%%%%%%%%%%%%%%%%%%%
%%%%%%%%%%%%%%%%%%%%%%%%%%%%%%%%%%%%%%%%%%%%%%%%%%%%%%%%%%%%%%%%%%%%
%%%%%%%%%%%%%%%%%%%%%%%%%%%%%%%%%%%%%%%%%%%%%%%%%%%%%%%%%%%%%%%%%%%%
%%%%%%%%%%%%%%%%%%%%%%%%%%%%%%%%%%%%%%%%%%%%%%%%%%%%%%%%%%%%%%%%%%%%
%%%%%%%%%%%%%%%%%%%%%%%%%%%%%%%%%%%%%%%%%%%%%%%%%%%%%%%%%%%%%%%%%%%%
%%%%%%%%%%%%%%%%%%%%%%%%%%%%%%%%%%%%%%%%%%%%%%%%%%%%%%%%%%%%%%%%%%%%

Let $\C[[x]]$ be the ring of formal power series with complex
coefficients and $\C\{x\}$ the ring of germs of analytic functions
at zero. A divergent formal power series $\hat f=\sum_{n\geq 0}f_n
x^{n+1}\in x\C[[x]]\smallsetminus x\C\{x\}$ is said to be a
generic Gevrey series if it is solution of a differential equation
of the form \begin{equation}\label{eq:equadiff}
a_n(x)\l(x^2\partial\r)^{n}y(x)+a_{n-1}(x)\l(x^2\partial\r)^{n-1}y(x)+
\dots+a_0(x)y(x)=g(x)\,, \end{equation} where $\partial=\frac
d{dx}$ and $a_0(x),\dots,a_n(x),g(x)\in\C\{x\}$, with
$a_0(0)a_n(0)\neq 0$. This implies that its formal Borel transform
$\cB(\hat f)=\sum_{n\geq 0}\frac{f_n}{n!}\xi^n\in\C\{\xi\}$ is a
germ of an analytic non entire function. The most important
example is the Euler series
$$
\hat E(x)=\sum_{n\geq 0}(-1)^nn!x^{n+1}\,,
$$
which is solution of the differential equation $x^2\partial
y+y=x$. A generic Gevrey series has the following properties:
$\cB(\hat f)$ can be analytically continued along almost all
direction $d\in(-\pi,\pi]$ and the Laplace integral along the half
line $e^{id}{\mathbb R}^+$:
$$
\cS^d(\hat f)=\int_0^{e^{id}\infty}\cB(\hat
f)(\xi)e^{-\xi/x}d\xi\,,
$$
called sum of $\hat f$ in the direction $d$, represents a
convergent solution of \eqref{eq:equadiff}, analytic on a
convenient sector and asymptotic to $\hat f$ at zero: this is the
first result of the well-known theory of summation of divergent
series (\cf
\cite{raminet,malgrangeexpositiones,lodaygazette,lodayexpositiones}).
\par
In the last fifteen years analogous summation theories for
$q$-difference equations have been developed (\cf
\cite{zhangfourier,MarotteZhang, zhanggroningen,RZ,DeSoleKac}). This
last sentence already shows one issue in the topic: there are many
$q$-summation theories in the literature and the relations among
them are not clear.
\par
Let us consider a $q$-deformation of the Euler series, namely:
$$
\hat E_q(x)=\sum_{n\geq 0}(-1)^n[n]_q^!x^{n+1}\,,
$$
where $[n]_q=1+q+\dots+q^{n-1}$ and
$[n]_q^!=[n]_q[n-1]_q\cdots[1]_q$. This series converges
coefficientwise to $\hat E(x)$ when $q\to1$ and is solution of the
$q$-difference equation
$$
x^2\dq y+y=x\,, \hbox{ with $q\in\C^*$ and $\ds\dq
y(x)=\frac{y(qx)-y(x)}{(q-1)x}$,}
$$
which is a discretization, in an obvious sense, of the so-called
Euler differential equation $x^2\partial y+y=x$. A first dichotomy
immediately appears: when $|q|<1$ the series $\hat E_q$ is a germ of
analytic function, converging for $|x|<|1-q|$, while for $|q|>1$ the
series $\hat E_q$ diverges. This is itself quite a curious fact,
that we have investigated in the present paper.
\par
As far as the divergent case $|q|>1$ is concerned, another
dichotomy immediately shows up: authors have been using two formal
Borel transforms, namely
$$
\cB_q(\hat f)=\sum_{n\geq 0}\frac{f_n}{[n]_q^!}\xi^n \hbox{ and }
B_q(\hat f)=\sum_{n\geq 0}\frac{f_n}{q^{n(n-1)/2}}\xi^n\,.
$$
Notice that we have $\cB_q(\hat E_q)=\frac 1{1+\xi}$ and $B_q(\hat
E_q)=\hat E_p(\xi)$, with $p=q^{-1}$. Each one of these formal Borel
transforms will be seen to naturally determine two summation
procedures, so that we end up with at least four summation
procedures: understanding the relations among them is a natural
question. Notice that from an arithmetic point of view, $\cB_q$ and
$B_q$ are deeply different (\cf \cite{andreannalsII}).

\begin{center}
$\ast\ast\ast$
\end{center}

The present paper is divided in two parts: in the first one we
consider the case $q\in(0,1)\subset{\mathbb R}$, while in the second
one we study different summation procedures under the assumption
$q\in(1,\infty)$. Let us make a few comments on these assumptions:
\begin{trivlist}
\item $\bullet$
We assume that the parameter
$q$ is real: this simplifies the exposition, although
it is not always completely necessary.
\item $\bullet$
Authors writing on
$q$-difference equations say sometimes that choosing $q$ smaller
or greater than one is only a matter of convention: as we explain
below, this is not true in the present situation, and the two
cases need to be investigated separately.
\end{trivlist}

\medskip
Let $q\in(0,1)\subset{\mathbb R}$. In this case $\hat E_q$ is the
Taylor expansion at $0$ of a meromorphic function $\cE_q$ on $\C$,
whose poles are a discrete subset of the negative real axis
${\mathbb R}^-$. In \S\ref{sec:qeuler-q<1} we prove the uniform
convergence of $\cE_q$ on the compacts of
$\C\smallsetminus{\mathbb R}^-$ to the analytic continuation $\cE$
of the Borel sum of $\hat E$ in the direction ${\mathbb R}^+$. The
proof of this result is based on the development of $\cE_q$ at
$\infty$, which is a $q$-deformation of the classical expansion of
$\cE$ at $\infty$ (\cite[page 261]{erd}):
$$
\cE(x)=(-\log x+\gamma)e^{\frac1x}+\sum_{n\ge 1}\frac{\sum_{1\le
k\le n}\frac1k}{n!}\l(\frac1x\r)^n\,,
$$
where  $\gamma=\lim_{n\to \infty}\l(\sum_{k=1}^n\frac 1k-\ln
(n)\r)$ is the Euler constant. In the same spirit, using Sauloy's
canonical solutions at $\infty$ of a fuchsian $q$-difference
operators \cite{Sfourier} and his result on their confluence when
$q\to 1$, we can prove the main theorem of the first part. Namely,
let $y(q,x)=\sum_{n\geq 0}y_n(q)x^{n+1}\in x\C[[x]]$ be a family
of formal power series, with $q\in(\eta,1]$, for some
$\eta\in(0,1)$. We suppose that the $y_n(q)$'s are continuous
functions of $q$ and that the family
$\phi(q,\xi)=\cB_qy(q,x)\in\C\{\xi\}$ is solution of a family of
equations over $\mathbb P^1_\C$, fuchsian and non resonant at
$\infty$\footnote{For a precise definition of a \emph{family of
equations over $\mathbb P^1_\C$, fuchsian and non resonant at
$\infty$}, \cf Assumption \ref{assumption} below.}. Then (\cf
Theorem \ref{thm:convergence1} below):

\begin{thmintro}\label{thmintro}
Let $d\in[0,2\pi)$ be such that $\phi(1,x)$ is holomorphic on a
domain containing the half line $e^{i d}{\mathbb R}^+$. Then  for
any $x\in V:=\{\vert\arg x-d\vert<\frac{\pi}{2}\}$ we have
$$
\lim_{q\to 1^-}y(q,x)=\cS^d(y(1,x))=
\int_0^{e^{id}\infty}\phi(1,\xi)e^{-\xi/x}d\xi\,,
$$
the convergence being uniform on the compacts of $V$.
\end{thmintro}

This result immediately implies two corollaries (\cf
\S\ref{subsec:applications} below). First of all, let
$y(x)=\sum_{n\geq 0}y_nx^{n+1}\in x\C[[x]]$ be a series such that
$\phi(\xi)=\sum_{n\geq 0}\frac{y_n}{n!}\xi^n$ is solution of a
fuchsian differential equation $\sum_{i=0}^\mu
A_i(\xi)(x\partial)^i \phi=0$ on $\mathbb P^1_\C$, non resonant at
$\infty$. One can construct a family of power series $\mathbf
y_q(x)$, with $q\in(0,1)$, such that $\cB_q\mathbf y_q(\xi)$ is
solution of $\sum_{i=0}^\mu A_i(\xi)(x\dq)^i\phi=0$ and $\mathbf
y_q(x)$ converges coefficientwise to $y(x)$ when $q\to 1^-$. Then:

\begin{corintro}
The family $\mathbf y_q(x)$ converges uniformly to the Borel sum
$\cS^d(y(\xi))$ of $y(x)$, when $q\to 1^-$, on the compacts of a
convenient sector $V=\{|\arg x-d|<\pi/2\}$.
\end{corintro}

A second corollary is about the sum of confluent
hypergeometric series. Let us consider $a,b\in\C$, with
$a-b\not\in{\mathbb Z}$, and the basic hypergeometric function:
$$
\Phi(a,b;q,x)=\sum_{n\geq 0}\frac{(q^a;q)_n(q^b;q)_n}{(q;q)_n}
\l(\frac{x}{1-q}\r)^n\,,
$$
where $(a;q)_0=1$ and $(a;q)_n=(1-a)(1-qa)\cdots(1-q^{n-1}a)$ for
any integer $n\geq 1$.

\begin{corintro}
The analytic function $\Phi(a,b;q,x)$ converges uniformly to the
Borel sum of the hypergeometric confluent series
$$
{}_2F_0(a,b;-;x)=\sum_{n\geq 0}\frac{(a)_n(b)_n}{n!}x^n\,,
$$
with $(a)_0=1$ and $(a)_n=a(a+1)\cdots(a+n)$ for any integer
$n\geq 1$, on the compacts of a convenient sector centered at $0$,
when $q\to 1^-$.
\end{corintro}

Finally, notice that the result on the confluence of $\cE_q$ can
be deduced from theorem \ref{thmintro}.

\medskip
The second part of the paper deals with the summation of divergent
$q$-series when $q\in(1,\infty)\subset{\mathbb R}$. Following the
scheme of the first part, we start our investigation studying four
summations of the series $\hat E_q(x)$. We consider the
$q$-exponential and the classical Theta function (here
$p=q^{-1}$):
$$
e_q(x)=\sum_{n\geq 0}\frac{x^n}{[n]_q^!} \,,\qquad
\theta_p(x)=\sum_{n\in{\mathbb Z}}p^{n(n-1)/2}x^n\,.
$$
We replace, in the classical Laplace integral, the exponential
function by $e_q$ or $\theta_p$.
By using an usual integral or a
discrete $q$-analogue\footnote{It's
precisely the so-called Jackson's integral, which is an infinite
sum that approximates the associated usual integral. The precise
definition is in Appendix \ref{sec:jacobsonintegral} below.},
denoted $\int_{\lambda p^{\mathbb Z}}f d_p\xi$, we
get four different {\it $q$-Borel sums of $\hat E _q(x)$}:
$$
\begin{array}{ll}
\ds \cE_q^d(x)=\frac{q-1}{\ln
q}\int_0^{e^{id}\infty}\frac1{(1+\xi)e_q(q\frac \xi x)}d\xi, &\ds
E_q^d(x)=\frac{q}{\ln q}\int_0^{e^{id}\infty}\frac{\hat E_p(\xi)}
    {\theta_p(q\frac \xi x)}d\xi;\\ \\
\ds\cE_q^{[\lambda]}(x)=\frac q{1-p}\int_{\la p^{\mathbb Z}}\frac1
    {(1+\frac \xi{1-p})e_q(q\frac \xi{(1-p)x})}\,{d_p\xi},
&\ds E_q^{[\lambda]}(x)=\frac q{1-p}
    \int_{\la p^{\mathbb Z}}\frac{\hat E_p(\xi)}{\theta_q(q\frac \xi
    x)}\,{d_p\xi},
\end{array}
$$
where $d\in(-\pi,\pi)$ and $\la\notin -p^{\mathbb Z}$. We prove
that $E^d_q(x)=\cE^d_q(x)$ on the sector $\arg(x)\in(-2\pi,2\pi)$
of the Riemann surface of the logarithm and that
$\cE_q^{[\la]}(x)=E_q^{[\la]}(x)$ for any
$x\in\C\smallsetminus(p-1)\la q^{\mathbb Z}$. Moreover we can
explicitly determine the functions $\cE^d_q(x)-\cE_q^{[\la]}(x)$
and $\cE_q^{[\la]}(x)-\cE_q^{[\mu]}(x)$ for
$x\in\C\smallsetminus{\mathbb R}^-$, in terms of the Theta
function. Finally, we establish the following relation between
$\cE^d_q(x)$ and $\cE_q^{[\la]}(x)$ (\cf Corollary
\ref{thm:moyenne} below):
$$
\cE_q^d(x)=\frac1{\ln q}\int_1^q\cE_q^{[\la]}(x)\frac{d\la}\la\,.
$$
\par
In an analogous way, for a formal power series $\hat f\in\C[[x]]$
such that $\cB_q\hat f$ is an analytic function with a
$q$-exponential growth of order one at $\infty$ we can define its
sums $\cS_q^d\hat f$, $\cS^{[\la]}_q\hat f$, $S_q^d\hat f$ and
$S^{[\la]}_q\hat f$.
Using some explicit results for the Tschakaloff series and
a $q$-convolution product adapted to the situation we can prove
the following (\cf Theorem \ref{thm:sums}):

\begin{thmintro}
Let $\hat f\in\C[[x]]$ be a generic $q$-Gevrey series. Then for
any $\la\in\C^\ast\smallsetminus\cup_{i=1}^n\mu_i q^{\mathbb Z}$,
for convenient $\mu_1,\dots,\mu_n\in\C^\ast$, and almost all
direction $d\in(-\pi,\pi)$, we have $\cS_q^d\hat f=S_q^d\hat f$
and $\cS^{[\la]}_q\hat f= S^{[\la]}_q\hat f$. Moreover:
$$
\ds\cS^d_q\hat f=\frac1{\ln q}\int_{e^{id}}^{qe^{id}}
\cS^{[\la]}_q\hat f\frac{d\la}\la\,.
$$
\end{thmintro}

\begin{center}
$\ast\ast\ast$
\end{center}

The theory of irregular singular $q$-difference equations is
nowadays relatively well understood. This paper deals with two of
the questions that are still without answer, namely:
\begin{trivlist}

\item
1. Thanks to the work of J. Sauloy \cite{Sfourier} we know how to
``uniformly approximate'' the global monodromy of a fuchsian
differential equation on $\mathbb P^1_\C$, in terms of the
Birkhoff matrices, which is a sort of $q$-monodromy, of a family
of $q$-difference equations deforming the given differential one.
Of course an analogous result is expected be true for the Stokes
phenomenon: actually the \emph{confluence} of the Stokes matrices
is studied for some functional equations linked to classical
special functions (\cf for instance \cite{zhanggroningen}).
The main theorem of the first part of this
article goes in the direction of a discrete deformation of the
Stokes phenomenon: differently from previous authors, we consider
a discrete convergent deformation of the divergent differential
datum. This approach is not really explored and the present result
surely does not exhaust its possible applications.

\item
2. In the second part of the paper we study the relations between
the different kind of $q$-Borel sums considered in the literature.
We prove the relations among them for a generic Gevrey series.
This is a first step towards the proof of a general result for a
divergent solution of a $q$-difference equations, having a Newton
polygon with more than one slope.

\end{trivlist}

\bigskip\noindent
{\bfseries Acknowledgement.} We would like to thank the organizers
of the S\'eminaire sur les \'equations aux $q$-diff\'erences, at
the University of Toulouse 3, and of the Special Session in
Differential Algebra, at the 2007 AMS Spring Eastern Sectional
Meeting, for giving us the possibility to expose the results
contained in this paper and their interest in our work. In
particular we thank Yves Andr\'e, Jean-Pierre Ramis and Jacques
Sauloy for their interest and encouragement.
\par
The first author would like to thanks the University of Lille 1,
and particularly the participants of the S\'eminaire de
th\'eorie de Galois diff\'erentielle, for their hospitality.
\par
The authors thanks the anonymous referee for his remarks and suggestions.

%%%%%%%%%%%%%%%%%%%%%%s%%%%%%%%%%%%%%%%%%%%%%%%%%%%%%%%%%%%%%%%%%%%%%
%%%%%%%%%%%%%%%%%%%%%%%%%%%%%%%%%%%%%%%%%%%%%%%%%%%%%%%%%%%%%%%%%%%%
%%%%%%%%%%%%%%%%%%%%%%%%%%%%%%%%%%%%%%%%%%%%%%%%%%%%%%%%%%%%%%%%%%%%
%%%%%%%%%%%%%%%%%%%%%%%%%%%%%%%%%%%%%%%%%%%%%%%%%%%%%%%%%%%%%%%%%%%%
%%%%%%%%%%%%%%%%%%%%%%%%%%%%%%%%%%%%%%%%%%%%%%%%%%%%%%%%%%%%%%%%%%%%
%%%%%%%%%%%%%%%%%%%%%%%%%%%%%%%%%%%%%%%%%%%%%%%%%%%%%%%%%%%%%%%%%%%%
\numberwithin{equation}{thm}

\section*{Part I. Convergent $q$-Borel and $q$-Laplace transform
and confluence: the case $q<1$}
\addcontentsline{toc}{section}{Part I. Convergent $q$-Borel and
$q$-Laplace transform and confluence: the case $q<1$}

%%%%%%%%%%%%%%%%%%%%%%%%%%%%%%%%%%%%%%%%%%%%%%%%%%%%%%%%%%%%%%%%%%%%
%%%%%%%%%%%%%%%%%%%%%%%%%%%%%%%%%%%%%%%%%%%%%%%%%%%%%%%%%%%%%%%%%%%%
%%%%%%%%%%%%%%%%%%%%%%%%%%%%%%%%%%%%%%%%%%%%%%%%%%%%%%%%%%%%%%%%%%%%
%%%%%%%%%%%%%%%%%%%%%%%%%%%%%%%%%%%%%%%%%%%%%%%%%%%%%%%%%%%%%%%%%%%%
%%%%%%%%%%%%%%%%%%%%%%%%%%%%%%%%%%%%%%%%%%%%%%%%%%%%%%%%%%%%%%%%%%%%
%%%%%%%%%%%%%%%%%%%%%%%%%%%%%%%%%%%%%%%%%%%%%%%%%%%%%%%%%%%%%%%%%%%%

We suppose that $q\in(0,1)\subset{\mathbb R}$ and we set $p=q^{-1}$.
\par
The first part of the paper is organized as follows. First of all
we study the properties of the $q$-deformation
$\cE_q(x)=\sum_{n\geq 0}(-1)^n[n]_q^!x^n$ of the Euler series:
namely we give two integral representations for $\cE_q(x)$, and
use them for proving that $\cE_q(x)$ converges uniformly to the
Borel sum of $\hat E(x)=\sum_{n\geq 0}(-1)^nn!x^n$ in the direction
${\mathbb R}^+$, uniformly on the compacts of a convenient sector.
Then we give an analogous result for general $q$-series, deforming
coefficientwise a Gevrey series of order 1. In appendix
\ref{sec:jacobsonintegral} we recall some general facts on the
Jackson integral, while in appendix \ref{sec:watsonformula} we
prove a degenerate $q$-Watson formula for Heine's series that we
need in \S\ref{sec:qeuler-q<1} for the proof of Proposition
\ref{prop:qEulerinfinity}.

%%%%%%%%%%%%%%%%%%%%%%%%%%%%%%%%%%%%%%%%%%%%%%%%%%%%%%%%%%%%%%%%%%%%
\section{Convergent $q$-Euler series}
\label{sec:qeuler-q<1}
%%%%%%%%%%%%%%%%%%%%%%%%%%%%%%%%%%%%%%%%%%%%%%%%%%%%%%%%%%%%%%%%%%%%

The series
$$
\cE_q(x)=\sum_{n\geq 0}(-1)^n[n]_q^!x^{n+1}\,,
$$
where $[n]_q=\frac{q^n-1}{q-1}$ and
$[n]_q^!=[n]_q[n-1]_q\cdots[1]_q$, represents a germ of analytic
function at $0$. If we consider the $q$-derivation:
$$
\dq y=\frac{y(qx)-y(x)}{(q-1)x}
$$
and observe that $\dq x^n=[n]_qx^{n-1}$ for any $n\in{\mathbb Z}$,
$n\geq 1$, then $\cE_q(x)$ verifies the functional equation :
$$
x^2\dq y+y=x\,,
$$
that can be rewritten in the form:
$$
y(x)=\frac{x}{x+1-q}y(qx)- \frac{(q-1)x}{x+1-q}\,.
$$
By substitution of $x$ by $q^nx$, we deduce that
$$
y(q^nx)=\frac{q^nx}{q^nx+1-q}y(q^{n+1}x)-
\frac{(q-1)q^nx}{q^nx+1-q}\,,
$$
which implies that $\cE_q(x)$ can be continued to an analytic
function on $\C\smallsetminus\{(q-1)q^n:n\in{\mathbb Z},\,n\leq
0\}$. The discrete spiral of poles $\{(q-1)q^n:n\in{\mathbb
Z},\,n\leq 0\}$ of $\cE_q$ turns out to be a spiral of simple
poles, as the following lemma shows:

\begin{lemma}
The analytic function $\cE_q$ admits the following expansion
\begin{equation}
\label{equa:Heine} \cE_q(x)=(1-q)\sum_{n\ge
0}\frac{(q^{n+1};q)_\infty}{1+\frac{1-q}{q^nx}}\,,
\end{equation}
where $(a;q)_\infty=\prod_{i=0}^\infty (1-q^ia)$.
\par
In particular for any $k\in{\mathbb Z}$, $k\leq 0$, the function
$\cE_q(x)$ has a simple pole at $(q-1)q^k$. The residue of the differential form
$\cE_q(x)dx$ at $(q-1)q^k$ is equal to
$$
Res_{x=(q-1)q^k}\cE_q(x) dx=-(1-q)^2q^k(q^{1-k};q)_\infty\,.
$$
\end{lemma}

We recall some standard notations for basic hypergeometric
functions
$$
\l\{
\begin{array}{l}
\ds {}_2\phi_1\l(a,b;c;q,x\r)=\sum_{n\geq 0}\frac{(a;q)_n(b;q)_n}{(c;q)_n(q;q)_n}x^n\,,\\
\hbox{where $(a;q)_0=1$ and $(a;q)_n=\prod_{k=1}^{n}(1-aq^{k-1}))$
for $1\le n\le \infty$, }
\end{array}\r.
$$
and the Heine's basic transformation (\cf \cite[\S 1.4]{GR}):
\begin{equation}\label{equa:heinetransformation} {}_2\phi_1\l(a,b;c;q,x\r)=
\frac{(a;q)_\infty(bx;q)_\infty}{(c;q;)_\infty(x;q)_\infty}
{}_2\phi_1\l(c/a,x;bx;q,a\r)\quad(\vert q\vert<1,\ \vert
a\vert<1)\,. \end{equation}

\begin{proof}
The lemma above is a straightforward application of
\eqref{equa:heinetransformation}, in fact:
\begin{equation}\label{equa:qeulerhypergeometric}
\cE_q(x)=x~{}_2\phi_1\l(q,q;0;q,-\frac x{1-q}\r)\,. \end{equation}
The calculation of the residues of $\cE_q(x)$ follows at once.
\end{proof}

%%%%%%%%%%%%%%%%%%%%%%%%%%%%%%%%%%%%%%%%%%%%%%%%%%%%%%%%%%%%%%%%%%%%
\subsection{Integral representation}
%%%%%%%%%%%%%%%%%%%%%%%%%%%%%%%%%%%%%%%%%%%%%%%%%%%%%%%%%%%%%%%%%%%%

Using the Jackson's integral (\cf Appendix
\ref{sec:jacobsonintegral} for the definition) we obtain the
following integral representation for $\cE_q$:

\begin{prop}
For any $x\in\C\smallsetminus\{(q-1)q^n:n\in{\mathbb Z},\,n\leq
0\}$, we have:
\begin{equation}\label{equation:q-Euler}
\cE_q(x) =\int_0^{\frac{x}{1-q}}\frac{(q(1-q)\frac
tx;q)_\infty}{t+1}d_qt =\int_{q^{{\mathbb
Z}}\frac{x}{1-q}}\frac{(q(1-q)\frac tx;q)_\infty}{t+1}d_qt.
\end{equation}
\end{prop}

\begin{proof}
Let us remark that $(q^{-k};q)_\infty=0$ for any $k\in{\mathbb
Z}$, $k\ge 0$. Then it follows from Remark
\ref{rmk:calculintegral}, that Formula \eqref{equa:Heine} is
equivalent to \eqref{equation:q-Euler}
\end{proof}

\begin{rmk}
A straightforward verification shows that the infinite product
$(q(1-q)x;q)_\infty$ represents a germ of analytic function at $0$
and that it verifies the equation
$$
y(px)=\l(1+(p-1)(-qx)\r)y(x)\,,
$$
or equivalently
$$
d_py(x)=-qy(x)\,.
$$
This implies that $(q(1-q)x;q)_\infty$ coincides with the analytic
function at $0$:
$$
e_p(-qx):=\sum_{n\geq 0}\frac{(-qx)^n}{[n]_p^!}\,,
$$
so that Equation \eqref{equation:q-Euler} takes the more familiar
shape:
$$
\cE_q(x)=\int_0^{\frac{x}{1-q}} \frac{e_p(-qt/x)}{t+1}d_qt\,,
$$
that so closely reminds the Euler integral:
$$
\cE(x)=\int_0^{+\infty}\frac{e^{-\frac tx}}{t+1}dt\,.
$$
The analytic function $\cE(x)$ can be continued to
$\C\smallsetminus{\mathbb R}^-$, it is asymptotic at zero to the Euler
series $\sum_{n\geq 0}(-1)^n n! x^{n+1}$ and is solution of the
differential equation $x^2y^\p+y=x$. In the following subsection
we are going to study the behavior of $\cE_q(x)$ with respect to
$\cE(x)$ when $q\to 1^-$.
\end{rmk}

%%%%%%%%%%%%%%%%%%%%%%%%%%%%%%%%%%%%%%%%%%%%%%%%%%%%%%%%%%%%%%%%%%%%
\subsection{Confluence}
\label{subsec:confluence}
%%%%%%%%%%%%%%%%%%%%%%%%%%%%%%%%%%%%%%%%%%%%%%%%%%%%%%%%%%%%%%%%%%%%

Let us denote by $\cE(x)$ the analytic continuation to
$\C\setminus (-\infty,0]$ of the Borel sum of $\hat E(x)$ in the
direction ${\mathbb R}^+$:
$$
\cE(x)=\int_0^{+\infty}\frac{e^{-\frac tx}}{t+1}dt,\quad \Re
x>0\,,
$$
and  by $\log x$ the analytic continuation to
$\C\setminus(-\infty,0]$ of $\log x$.

\begin{thm}\label{thm:convergence}
If $q\to 1^-$, the analytic continuation of $\cE_q(x)$ converges
to $\cE(x)$ for any $x\in\C\smallsetminus(-\infty,0]$ and the
convergence is uniform on the compacts of
$\C\smallsetminus(-\infty,0]$.
\end{thm}

The proof of the theorem above relies on the following result (\cf
\S \ref{subsec:ProofProp} below for the proof):

\begin{prop}
\label{prop:qEulerinfinity} The following identity holds, for any
$x\in\C^*\setminus q^{-{\mathbb N}}$:
\begin{equation}\label{equa:qEulerinfinity}
\sum_{n\ge 0}(q;q)_nx^{n+1}
=-\l(-qx\frac{\theta^\prime(-qx)}{\theta(-qx)}+1+A(q)\r) \l(\frac
qx;q\r)_\infty +\sum_{n\ge
1}\frac{a_{n}}{(q;q)_{n}}q^{n(n+1)/2}\l(-\frac1x\r)^{n}\,,
\end{equation}
where
$$
\theta(x)=\theta(q,x)=\sum_{n\in{\mathbb Z}}q^{n(n-1)/2}x^n\,,
$$
$$
\ds A(q)=\sum_{n\ge 0}\frac{q^{n+1}}{q^{n+1}-1}\,
$$
and
$$
\ds a_{n+1}=\sum_{k=0}^n\frac1{q^{k+1}-1},\quad n\ge 0\,.
$$
\end{prop}

Our strategy for the proof of Theorem \ref{thm:convergence} is
based on the fact that \eqref{equa:qEulerinfinity} is a
``deformation'' of the following classical formula:
\begin{equation}\label{equa:Eulerinfinity}
\cE(x)=(-\log x+\gamma)e^{\frac1x}+\sum_{n\ge 1}\frac{\sum_{1\le
k\le n}\frac1k}{n!}\l(\frac1x\r)^n\,,
\end{equation}
where  $\gamma$ is the Euler constant:
$$
\gamma=\lim_{n\to \infty}\l(\sum_{k=1}^n\frac 1k-\ln (n)\r)\,.
$$
In fact, taking the logarithmic derivative of the functional
equation $\theta(x)=x\theta(qx)$, one proves that the
meromorphic function $(q-1)z\frac{\theta^\p(-z)}{\theta(-z)}$
verifies the equation $y(qx)-y(x)=q-1$ or equivalently $\dq
y(x)=\frac{1}{x}$, therefore it ``deforms'' the logarithm. On the
other hand we have:
$$
(q-1)A(q)=\sum_{n\geq 0}\frac{q^{n+1}}{[n+1]_q}\,,
$$
whose link to the Euler constant is intuitive. The proof of
Theorem \ref{thm:convergence} is a formalization of these ideas.

%\begin{parag}
%{\bfseries
\subsubsection*{Proof of Theorem \ref{thm:convergence}}

If we perform the variable change $x\to \frac{x}{q-1}$ in
\eqref{equa:Eulerinfinity} and remember that
$$
e_p\l(q/x\r)=\l(\frac{-q(1-q)}x;q\r)_\infty,
\quad \ell_q(x):=-x\frac{\theta'(-x)}{\theta(-x)}\,,
$$
then we
obtain the expression
\begin{equation}\label{equa:Eulerinfinitybis}
\cE_q(x) =(q-1)\,\l[\ell_q\l(\frac{qx}{q-1}\r)
+1+A(q)\r]\,e_p\l(\frac qx\r) +\sum_{n\geq
1}\frac{(q-1)a_n}{[n]_q^!}q^{n(n-1)/2}\l(\frac qx\r)^n\,,
\end{equation}
that we are going to analyze term by term.

\par
First of all the constant $A(q)$ can be expressed in terms of the
logarithmic derivative
$\Psi_q(x)=\frac{\Gamma_q^\p(x)}{\Gamma_q(x)}$ (see
\eqref{equa:PsiOmega}, where $\Omega(q)=A(q)$):
\begin{equation}\label{equa:A(q)Psi}
A(q)=-\frac1{\ln q}\l(\Psi_q(1)+\ln(1-q)\r).
\end{equation}

The following result says how the $q$-logarithm $\ell_q$ tend to
the usual logarithm.

\begin{lemma}\label{lemma:log}
Let $\epsilon\in(0,\pi)$ and consider the sector
$V_\epsilon=\{x\in\C^*: \vert\arg x\vert\le \pi-\epsilon\}$. Then
the following uniform estimate holds for any $(q,x)\in(0,1)\times
V_\epsilon$:
\begin{equation}\label{equa:qlogestimate}
\vert\ln q\,\ell_q(-\sqrt q\,x)+\log x\vert\le
\frac{4\pi\,e^{\frac{2\pi}{\ln
q}\,\epsilon}}{(1-e^{\frac{4\pi^2}{\ln q}})(1-e^{\frac{2\pi}{\ln
q}\,\epsilon})} \,.
\end{equation}
\end{lemma}

\begin{proof}
The lemma is a consequence of the following classical functional
relation for $\theta(x)$ (\cf \cite[\S 21.51, p. 475]{WW}, where
$\vartheta_3(z\vert\tau)=\theta(\sqrt{q}\,e^{2\pi iz})$ with $q=e^{2\pi i\tau}$):
\begin{equation}\label{equa:Modular}
\theta(\sqrt q\,
x)=\sqrt{\frac{2\pi}{\ln(1/q)}}\,e^{-\frac{\log^2x}{2\ln
q}}\theta^*(\sqrt{q^*}\,x^*),
\end{equation}
where we write
$$
x^*=e^{-\frac{2\pi i}{\ln q}\log x},\quad q^*=e^{\frac{4\pi^2}{\ln
q}}
$$
and denote by $\theta^*$ the Theta function obtained by replacing
$q$ by $q^*$. Indeed, if we take the logarithmic derivative w.r.t.
the variable $x$ in \eqref{equa:Modular} and observe that $\ln
q\,xdx^*=-2\pi ix^*dx$, then we obtain the following expression:
\begin{equation}\label{equa:Modular1}
\ln q\,\ell_q(-\sqrt q\,x)+\log x=-2\pi
i\ell_{q^*}(-\sqrt{q^*}\,x^*),
\end{equation}
so that we only need to examine $\ell_{q^*}(-\sqrt{q^*}\,x^*)$.
The key point of the proof is the fact that $q^*\to 0^+$
when $q\to 1^-$.
\par
For $\epsilon\in(0,\pi)$ we set:
$$
r_\epsilon=e^{\frac{2\pi}{\ln q}\epsilon}\in(\sqrt{q^*}\,,1),\quad
V_\epsilon^*=\{x\in\C:\frac{\sqrt{q^*}}{r_\epsilon}\,\le \vert
x\vert\le \frac{r_\epsilon}{\sqrt{q^*}}\,\}.
$$
It's obvious that
for any $x\in V_\epsilon$, we have $x^*\in V_\epsilon^*$, so that
\begin{equation}\label{equa:Module}
q^*<\frac{q^*}{r_\epsilon}\le \vert\sqrt{q^*}\,x^*\vert\le
r_\epsilon<1.
\end{equation}
On the other hand, the following identity, consequence of the Jacobi triple product,
$$
X\,\frac{{\theta^*}^\prime(X)}{\theta^*(X)}=\sum_{n\ge
0}\l(\frac{{q^*}^nX}{1+{q^*}^nX}-\frac{{q^*}^{n+1}}{X+{q^*}^{n+1}}\r)\,,
$$
combined with the inequality:
$$
\vert 1+X{q^*}^n\vert\ge 1-\vert X\vert,\quad
\vert X+{q^*}^{n+1}\vert\ge \vert X\vert-q^*\,,
\hbox{~for $q^*<\vert X\vert<1$,}
$$
implies that:
\begin{equation}\label{equa:qlogestimate1}
\sup_{\frac{q^*}{r_\epsilon}\le \vert X\vert\le
r_\epsilon}\l\vert\ell_{q^*}(-X)\r\vert
\le\frac{r_\epsilon}{1-q^*}\,\frac1{1-r_\epsilon}+
\frac{q^*}{1-q^*}\,\frac1{\frac{q^*}{r_\epsilon}-q^*}
\le\frac{2r_\epsilon}{(1-r_\epsilon)(1-q^\epsilon)}\,.
\end{equation}
We get \eqref{equa:qlogestimate} and hence Lemma \ref{lemma:log}
by combing \eqref{equa:Modular1} and \eqref{equa:qlogestimate1}.
\end{proof}

\begin{proof}[End of the proof of Theorem \ref{thm:convergence}]
By replacing $x$ by $\sqrt qx/(1-q)$ in Lemma \ref{lemma:log}, we
have:
$$
\ell_q\l(\frac {qx}{q-1}\r)=-\frac1{\ln q}\,
\l[\log x+\ln\l(\frac{\sqrt q}{1-q}\r)+O\l(e^{2\pi\epsilon/\ln q}\r)\r],
$$
so that we obtain, by \eqref{equa:A(q)Psi},
\begin{equation}\label{equa:ellA}
\ell_q\l(\frac {qx}{q-1}\r)+1+A(q)=-\frac1{\ln q}\,\l[\log
x+\Psi_q(1)+\frac{\ln q}2+O(e^{2\pi\epsilon/\ln q})\r]\,,
\end{equation}
where $\Psi_q$ denotes the logarithmic derivative of $\Gamma_q$.
As $q\to 1^-$, the function $\Gamma_q(x)$ converges uniformly to
$\Gamma(x)$ on any compact of $\C\setminus(-{\mathbb N})$ (\cf
\cite{zhangqGamma}), so $\Psi_q(x)$ converges to the logarithmic
derivative $\Psi(x)$ of the $\Gamma$ function. From the classical
relation $\Psi(1)=-\gamma$, one deduces that
$\Psi_q(1)=-\gamma+o(1)$. In other words,  \eqref{equa:ellA}
implies the following estimate:
\begin{equation}\label{equa:ellA1}
(q-1)\,\l[\ell_q\l(\frac {qx}{q-1}\r)+1+A(q) \r]=-\log
x+\gamma+o(1)\,
\end{equation}
where $o(1)$ denotes a quantity converging to $0$ as $q\to 1^-$,
uniformly on any compact of $\C\setminus(-\infty,0]$.

Notice that the exponential function $e^{\frac{1}{x}}$ is the
uniform limit on any domain $\{\vert x\vert> R>0\}$ of the
$p$-exponential $e_p(q/x)$, for $$ e_p(\frac{1}{x})=\sum_{n\ge
0}\frac{q^{n(n-1)/2}}{[n]_q^!}\,\l(\frac qx\r)^n\,.
$$
In the same time, again the dominated convergence Theorem implies
that, as $q\to 1^-$,
$$
\sum_{n\ge 1}\frac{(q-1)a_nq^{n(n-1)/2}}{[n]_q^!}\l(\frac qx\r)^n\to
\sum_{n\ge 1}\frac{\sum_{k=1}^n\frac1k}{n!}\l(\frac1x\r)^n\,,
$$
uniformly for $\vert x\vert> R>0$.
We
conclude combining \eqref{equa:Eulerinfinitybis} with
\eqref{equa:ellA1}.
\end{proof}
%\end{parag}

%%%%%%%%%%%%%%%%%%%%%%%%%%%%%%%%%%%%%%%%%%%%%%%%%%%%%%%%%%%%%%%%%%%%%%%%%%%%
%%%%%%%%%%%%%%%%%%%%%%%%%%%%%%%%%%%%%%%%%%%%%%%%%%%%%%%%%%%%%%%%%%%%%%%%%%%%
%%%%%%%%%%%%%%%%%%%%%%%%%%%%%%%%%%%%%%%%%%%%%%%%%%%%%%%%%%%%%%%%%%%%%%%%%%%%
%%%%%%%%%%%%%%%%%%%%%%%%%%%%%%%%%%%%%%%%%%%%%%%%%%%%%%%%%%%%%%%%%%%%%%%%%%%%
%%%%%%%%%%%%%%%%%%%%%%%%%%%%%%%%%%%%%%%%%%%%%%%%%%%%%%%%%%%%%%%%%%%%%%%%%%%%

\section{Confluence of the \emph{convergent} $q$-analogue of Borel-Laplace summation}

Let $q$ be a real number in the open interval $(0,1)$. We want to
generalize, under convenient reasonable assumptions, the results
of the previous section.

\subsection{Definition of the convergent $q$-Borel and
$q$-Laplace transform}

\begin{defn}
Let $\C\{x\}$ be the ring of the germs of analytic functions in
the neighborhood of $x=0$.
\begin{enumerate}
\item
We call \emph{(convergent) $q$-Borel transform} the map ${\cal
B}_q$ given by:
$$
{\cal B}_q\ :\quad x\C\l\{x\r\}\to \C\{\xi\},\ \sum_{n\ge 0}a_n
x^{n+1}\mapsto \sum_{n\ge 0}\frac{a_n}{[n]_q^!}\xi^n.
$$
\item
The \emph{(convergent) $q$-Laplace transform ${\cal L}_q$} is
defined by
$$
{\cal L}_q={\cal B}_q^{-1}\quad :\quad \C\{\xi\}\to x\C\l\{x\r\},\
\sum_{n\ge 0}a_n\xi^n\mapsto \sum_{n\ge 0}a_n[n]_q^!x^{n+1}.
$$
\end{enumerate}
\end{defn}

\begin{rmk}
Notice that the $q$-Euler series $\cE_q(x)$, considered in the
previous section, converges for $|x|<1-q$. Therefore a function
$f(x)$ is analytic on an open disc $\{\vert x\vert<R\}$, for some
$R\in(0,\infty)$, if and only if its $q$-Borel transform ${\cal
B}_qf(\xi)$ is analytic for  $\vert \xi\vert<R/(1-q)$.
\par
Calling $\cB_q$ and $\cL_q$ $q$-Borel and $q$-Laplace transform is
somehow an abuse of language: they don't transform convergent
series in divergent series and \emph{vice versa}. Nevertheless
they have interesting properties and we will show that they play a
role in the understanding of the confluence in the irregular case.
In fact, when $q\to1$, they tend coefficientwise to the usual
Borel and Laplace transforms, that we will denote $\cB_1$ and
$\cL_1$ respectively.
\end{rmk}

An important property of $\cB_q$ and $\cL_q$ is that they can be
expressed both as continuous and discrete integrals:

\begin{prop}\label{prop:fg}
Let $f\in x\C\{x\}$ and $\phi\in\C\{\xi\}$ such that ${\cal
B}_qf=\phi$. Then:
$$
\begin{array}{c}
\ds\phi(\xi)=\frac{1}{2\pi i}\int_{\vert x\vert=R} \frac
{f(x)}{((1-q)\frac{\xi}{x};q)_\infty }\frac{dx}{x^2}
=\frac{1}{2\pi i}\int_{\vert x\vert=R}f(x)e_q\l(\xi/x\r)\frac{dx}{x^2},\\ \\
\ds f(x)=\frac{-1}{2\pi i}\int_{\vert
\xi\vert=\rho}\phi(\xi)\cE_q\l(-\frac{x}{\xi}\r)d\xi,
\end{array}
$$
where the radius $R$ and $\rho$ are assumed to be chosen
sufficiently small.
\end{prop}

\begin{proof}
The first equality is a consequence of the identity
$$
\frac{1}{(x;q)_\infty}=\sum_{n\ge 0}\frac{1}{(q;q)_n}x^n\,
$$
and of the residue theorem.
Taking into account \eqref{equa:Heine}, the second equality is an
application of the residue theorem.
\end{proof}

\begin{cor}
Let $f$ and $\phi$ be as in Proposition \ref{prop:fg}. Then:
$$
\phi(\xi)=(q;q)_\infty\int_0^\xi f((1-q)x)
\frac{(\frac{qx}{\xi};q)_\infty}{\theta'(-\frac{x}{\xi})}d_qx,
\quad
$$
\begin{equation}\label{equa:fintq}
f(x)
=\int_0^{\frac{x}{1-q}}\l(\frac{(1-q)q\xi}{x};q\r)_\infty\phi(\xi)d_q\xi
=\int_0^{\frac{x}{1-q}}e_q\l(qx/\xi\r)^{-1}\phi(\xi)d_q\xi\,.
\end{equation}
\end{cor}

\begin{rmk}
Notice that Formula \eqref{equa:fintq} generalizes
\eqref{equation:q-Euler} and can be obtained directly by identifiying the coefficients. We give an alternative proof below.
\end{rmk}

\begin{proof}
Taking the derivative with respect to $x$ of the functional equation
$$
\theta(q^nx)=x^{-n}q^{-n(n-1)/2}\theta(x)\,
$$
and setting $x=-1$,
we obtain
$$
\theta'(-q^n)=(-1)^nq^{-n(n+1)/2}\theta'(-1)=(-1)^nq^{-n(n+1)/2}(q;q)_\infty^3\,.
$$
Again the residues formula and Equation \eqref{equa:Heine} imply
that
$$
\phi(\xi)=\frac{(1-q)\xi}{(q;q)_\infty}\sum_{n\ge
0}(-1)^n\frac{f((1-q)\xi q^n)}{(q;q)_n}q^{n(n+3)/2}
$$
and that
$$
f(x)=x\sum_{n\ge 0}q^n(q^{n+1};q)_\infty
\phi\l(\frac{q^nx}{1-q}\r)\,.
$$
This ends the proof.
\end{proof}

%%%%%%%%%%%%%%%%%%%%%%%%%%%%%%%%%%%%%%%%%%%%%%%%%%%%%%%%%%%
\subsection{Main result}
%%%%%%%%%%%%%%%%%%%%%%%%%%%%%%%%%%%%%%%%%%%%%%%%%%%%%%%%%%%

The formulas above suggest the convergence of the $q$-Laplace
transform ${\cal L}_q\phi$ to the classical Laplace transform
${\cal L}^d\phi$ (in the direction $d\in (0,2\pi)$):
\begin{equation}\label{equa:Laplace}
{\cal L}^d\phi(x)=\int_0^{\infty
e^{id}}\phi(\xi)e^{-\frac{\xi}{x}}d\xi\,,
\end{equation}
where $\phi$ is supposed to be holomorphic in a neighborhood of
$\xi =0$ and to be analytically continued in an open sector
$\{\vert\arg \xi-d\vert<\epsilon\}$ with at most an exponential
growth at infinity.

\begin{thm}\label{thm:convergence1}
Let $y(q,x)=\sum_{n\geq 0}y_n(q)x^{n+1}\in x\C[[x]]$ be a family
of formal power series, with $q\in(\eta,1]$, for some
$\eta\in(0,1)$. We suppose that the $y_n(q)$'s are continuous
functions of $q$ and that the family
$\phi(q,\xi)=\cB_qy(q,x)\in\C\{\xi\}$ is solution of a family of
equations over $\mathbb P^1_\C$, fuchsian and non resonant at
$\infty$, in the sense of Assumption \ref{assumption} below.
\par
Let $d\in[0,2\pi)$ be such that $\phi(1,x)$ is holomorphic on a
domain containing the half line $[0,e^{i d}\infty)$. Then  for any
$x\in V:=\{\vert\arg x-d\vert<\frac{\pi}{2}\}$ we have
$$
\lim_{q\to 1^-}y(q,x)={\cal L}^d\phi(1,\xi)\,,
$$
the convergence being uniform on any compact of $V$.
\end{thm}

Notice that $y(q,x)={\cal L}_q\phi(q,\xi)$, so that the result
above is actually a result about the confluence of $q$-summation.
Moreover $\phi(q,\xi)$ is meromorphic over $\C^\ast$ and its poles
are contained in a finite set of lines passing through the origin.
Also for $\phi(1,\xi)$ there are only a finite numbers of direction
$d$ that are forbidden: the anti-Stokes directions.

\begin{assumption}\label{assumption}
We suppose that:
\begin{enumerate}

\item
The series $\phi(1,\xi)$ is solution of a differential equation
$\cN_1\phi(1,\xi)=\sum_{i=0}^\mu A_i(1,\xi)\delta^i\phi(1,\xi)=0$,
where $\delta=\xi\frac{d}{d\xi}$, $A_i(1,\xi)\in\C[\xi]$, and the
operator $\cN_1$ is fuchsian at $0$ and $\infty$. Moreover we
suppose that the exponents of $\cN_1\phi(1,\xi)=0$ at $\infty$ are
non resonant.

\item
The series $\phi(q,\xi)=\cB_qy(q,x)$, $q\in(0,1)$, are solutions
of a linear $q$-difference operator
$\cN_q\phi(q,\xi)=\sum_{i=0}^\mu A_i(q,\xi)\delta_q^i
\phi(q,\xi)=0$, where $\delta_q=\xi\dq$, $A_i(q,\xi)\in\C[\xi]$,
and $\cN_q$ is fuchsian at $0$ and $\infty$\footnote{This means
that the only non vertical slope of the Newton-Ramis polygon, \ie
of the convex envelope of the set
$$
\{(i,j):~A_i(q,x)\neq 0,~\ord_{x=0}A_i(q,x)\leq j\leq \deg_x A_i(q,x)\}\,,
$$
are horizontal.}.

\item
The Newton-Ramis polygons of $\cN_q$ coincide for any
$q\in(\eta,1]$, and the coefficients $A_i(q,\xi)$ tends uniformly
to $A_i(1,\xi)$ when $q\to 1$, on any compact of $\mathbb P^1_\C$.
This implies in particular that for $q$ sufficiently closed to
$1$, the exponents of $\cN_q$ at $\infty$ are non resonant.

\item
For any $q$ sufficiently closed to $1$ there exists a constant
gauge transformation $C(q)\in Gl_\mu(\C)$ such that the constant
term at $\infty$ of the matrix
$$
C(q)^{-1}\l(\begin{array}{c|ccc}
0 & 1 &&0\\
\vdots &  &\ddots \\
0&0 & & 1\\
\hline -\frac{A_0(q,x)}{A_\mu(q,x)}&
-\frac{A_1(q,x)}{A_\mu(q,x)}&\dots&-\frac{A_{\mu-1}(q,x)}{A_\mu(q,x)}
\end{array}\r)
C(q)
$$
is in the Jordan normal form. We suppose that for $q\in(\eta,1]$
the entries of the matrix $C(q)$ are continuous functions of $q$
and that the form of the Jordan blocks is independent of $q$.

\end{enumerate}

\end{assumption}

%%%%%%%%%%%%%%%%%%%%%%%%%%%%%%%%%%%%%%%%%%%%%%%%%%%%%%%%%%%%%%%%%%%%%%%%%%%%
%%%%%%%%%%%%%%%%%%%%%%%%%%%%%%%%%%%%%%%%%%%%%%%%%%%%%%%%%%%%%%%%%%%%%%%%%%%%
%%%%%%%%%%%%%%%%%%%%%%%%%%%%%%%%%%%%%%%%%%%%%%%%%%%%%%%%%%%%%%%%%%%%%%%%%%%%
\subsection{Applications}
\label{subsec:applications}
%%%%%%%%%%%%%%%%%%%%%%%%%%%%%%%%%%%%%%%%%%%%%%%%%%%%%%%%%%%%%%%%%%%%%%%%%%%%
%%%%%%%%%%%%%%%%%%%%%%%%%%%%%%%%%%%%%%%%%%%%%%%%%%%%%%%%%%%%%%%%%%%%%%%%%%%%

Notice that the assumptions of Theorem \ref{thm:convergence1} are
verified in the following two natural situations.

\begin{parag}
[\bfseries ``Constant coefficient deformation'' of a differential
equation] For a linear differential equation $\sum_{i=0}^\mu
A_i(x)\delta^i y=0$, a possible trivial deformation is given by
$\sum_{i=0}^\mu A_i(x)\delta_q^i y=0$. One verifies that if
$\sum_{i=0}^\mu A_i(\xi)\delta^i y=0$ satisfies the first point of
Assumption \ref{assumption}, then $\sum_{i=0}^\mu
A_i(\xi)\delta_q^i y=0$ verifies automatically the next three
assumptions, provided that $1-q$ is small enough. Therefore we
have:

\begin{cor}
Let $y(x)=\sum_{n\geq 0}y_nx^{n+1}\in x\C[[x]]$ be a Gevrey series
of order one such that $\phi(\xi)=\cB_1 y(x)$ is solution of a
fuchsian differential equation $\sum_{i=0}^\mu A_i(x)\delta^i
\phi=0$ on $\mathbb P^1_\C$, non resonant at $\infty$. Consider a
family of power series $\mathbf y_q(x)$, with $q\in(0,1)$, such
that $\cB_q(\mathbf y_q)(\xi)$ is solution of $\sum_{i=0}^\mu
A_i(\xi)\delta_q^i\phi=0$ and $\mathbf y_q(x)$ converges
coefficientwise to $y(x)$ when $q\to 1^-$.
\par
Then the family $\mathbf y_q(x)$ converges uniformly to the Borel
sum of $y(x)$, when $q\to 1^-$, on the compacts of a convenient
sector $V=\{|\arg x-d|<\pi/2\}$.
\end{cor}

\end{parag}

\begin{parag}
[\bfseries Confluent hypergeometric case] Take $\phi(q,\xi)$ to be
the basic hypergeometric series:
$$
\phi(q,\xi)=\l\{\begin{array}{ll} \ds
{}_2\Phi_1(q^a,q^b;q;q,x)=\sum_{n\geq
0}\frac{(q^a;q)_n(q^b;q)_n}{(q;q)_n(q;q)_n}x^n\,,
& \hbox{ for $q\in(0,1)$};\\ \\
\ds {}_2F_1(a,b;1;x)=\sum_{n\geq 0}\frac{(a)_n(b)_n}{n!n!}x^n\,,&
\hbox{ if $q=1$;}
\end{array}\r.
$$
where $a,b\in\C$, with $a-b\not\in{\mathbb Z}$. Then Theorem
\ref{thm:convergence1} says that:

\begin{cor}
The basic hypergeometric analytic function
$$
\sum_{n\geq 0}\frac{(q^a;q)_n(q^b;q)_n}{(q;q)_n}
\l(\frac{x}{1-q}\r)^n
$$
converges uniformly to the Borel sum of the hypergeometric
confluent series
$$
{}_2F_0(a,b;-;x)=\sum_{n\geq 0}\frac{(a)_n(b)_n}{n!}x^n
$$
on the compacts of a convenient sector centered at $0$, when $q\to
1^-$.
\end{cor}

Of course the results above can be generalized. In fact, for any
$\ell\geq 2$, and any \emph{generic} choice of the parameters
$a_1,\dots,a_\ell,b_1,\dots,b_{\ell-2}\in\C$, the analytic basic
hypergeometric function
$$
\sum_{n\geq 0}\frac{(q^{a_1};q)_n\cdots(q^{a_\ell};q)_n}
{(q^{b_1};q)_n\cdots(q^{b_{\ell-2}};q)_n(q;q)_n}\l(\frac{x}{1-q}\r)^n
$$
converges uniformly to the Borel sum of the hypergeometric
confluent series
$$
{}_\ell F_{\ell-2}(a_1,\dots,a_\ell;b_1,\dots,b_{\ell-2};x)
=\sum_{n\geq 0}\frac{(a_1)_n\cdots(a_\ell)_n}{(b)_1\cdots
(b)_{\ell-2}n!}x^n
$$
on the compacts of a convenient sector centered at $0$, when $q\to
1^-$.
\end{parag}

%%%%%%%%%%%%%%%%%%%%%%%%%%%%%%%%%%%%%%%%%%%%%%%%%%%%%%%%%%%%%%%%%%%%%%%%%%%%
%%%%%%%%%%%%%%%%%%%%%%%%%%%%%%%%%%%%%%%%%%%%%%%%%%%%%%%%%%%%%%%%%%%%%%%%%%%%
%%%%%%%%%%%%%%%%%%%%%%%%%%%%%%%%%%%%%%%%%%%%%%%%%%%%%%%%%%%%%%%%%%%%%%%%%%%%
\subsection{Proof of Theorem \ref{thm:convergence1}}
%%%%%%%%%%%%%%%%%%%%%%%%%%%%%%%%%%%%%%%%%%%%%%%%%%%%%%%%%%%%%%%%%%%%%%%%%%%%
%%%%%%%%%%%%%%%%%%%%%%%%%%%%%%%%%%%%%%%%%%%%%%%%%%%%%%%%%%%%%%%%%%%%%%%%%%%%

%Notice that a variable change of the form $\xi\mapsto \la\xi$ allows to suppose that
%$d$ coincides with positive real semiaxis.
We know that our germs  $\phi(q,\xi)$, $q\in(\eta,1]$, of analytic
functions at $0$ admit an analytic continuation along $d$.
Moreover, for $q<1$, the functions $\phi(q,\xi)$ are actually
meromorphic over $\C$, which means that they are linear
combination of a basis of solutions of $\cN_qy=0$ at $\infty$. The
main point of the proof is the careful choice of such a basis,
that will allow us to prove that $\phi(q,x)$ converges uniformly
to $\phi(1,\xi)$ on an \emph{infinite} sector containing the
direction $d$. Of course this ends the proof since Equations
\eqref{equa:fintq} and \eqref{equa:Laplace} imply that for any
$x\in V$ we have:
\begin{equation}
\begin{array}{rcl}
\ds\lim_{q\to 1} y(q,x) &=&\ds\lim_{q\to1}
    \int_0^{\frac{x}{1-q}}\l(\frac{(1-q)q\xi}{x};q\r)_\infty\phi(q,\xi)d_q\xi\\ \\
&=&\ds\int_0^{\infty e^{i\arg(x)}}
    \lim_{q\to 1}\l(e_q\l(\frac{q\xi} x\r)^{-1}\phi(q,\xi)\r)d\xi\\ \\
&=&\ds\int_0^{\infty e^{i\arg(x)}}\phi(1,\xi)e^{-\frac\xi
x}d\xi\,.
\end{array}
\end{equation}

The theorem results of the combination of two lemmas. First of all
let us prove the uniform convergence around zero:

\begin{lemma}\label{lemma:uniformconvergenceat0}
The family $\phi(q,\xi)$ converges uniformly to $\phi(1,x)$, when
$q\to1$, on a closed disk centered at $0$.
\end{lemma}

\begin{proof}%[Proof of Lemma \ref{lemma:uniformconvergenceat0}]
Let us write $\phi(q,\xi)=\sum_{n\geq 0}\phi_n(q)x^n$ for any
$q\in(\eta,1]$. Then there exists $N>0$ such that for any $n>N$
the coefficients $\phi_n(q)$ verify a well defined recursive
relation whose coefficients do not degenerate\footnote{Notice that
we have not assumed that $0$ is a regular singular point with non
resonant exponents at $0$.}. A direct estimates of the recursive
relation allows to conclude that $|\phi_n(q)|\leq C^n$ for a
convenient real constant $C$, any $n>N$ and any $q\in(\eta,1]$.
This estimate, together with the fact that $\phi_n(q)$ is a
continuous function of $q$, implies the uniform convergence on a
convenient closed disk centered at $0$ (\cf for instance the
estimates in \cite[Lemma
1.2.6]{Spoligono}\nocite{sauloyfiltration}).
\end{proof}

The last assumption in \ref{assumption} implies that $\phi(q,x)$
is a linear combination, whose coefficients are entries of the
matrix $C(q)$, of the canonical solutions at $\infty$, constructed
in \cite[\S1]{Sfourier}, using a $q$-analog of the Frobenius
method. As noticed in \cite[\S3]{Sfourier}, the uniform part of
the canonical solution at $\infty$ converges uniformly on any
compact of $\mathbb P^1_\C\smallsetminus\{0\}$ where it is
analytic, to the uniform part of the solutions of $\cN_1y=0$
constructed with the classical Frobenius methods, once that the
gauge transformation $C(1)$ has been applied to the companion
matrix. Since the entries of $C(q)$ converges to the entries of
$C(1)$ by assumption, to obtain the uniform convergence in a
neighborhood of $\infty$ it is enough to control the convergence
of the so called \emph{log-car} matrix\footnote{The terminology
comes from the juxtaposition of the terms ``logarithm'' and
``character'', meaning the solution matrix of a constant
coefficient differential (resp. $q$-difference) system is obtained
by a combinatoric procedure from the logarithm (resp.
$q$-logarithm) and a family of characters (resp. $q$-characters).
A solution in a regular singular point, whose exponents are non
resonant, is given by the product of an analytic matrix,
called uniform part, by the ``log-car'' matrix.
\par
We are choosing here as a $q$-logarithm the logarithmic derivative
of the Jacobi $\theta$ function and as $q$-characters convenient
quotient of the $\theta$ functions. For more details in the
$q$-difference setting \cf \cite{Sfourier}.}. Let $\zeta=1/\xi$.
The uniform convergence of
$(1-q)\zeta\theta^\p(q,\zeta)/\theta(q,\zeta)$ over the infinite
sector $\{|\arg(\zeta)|<\pi-\veps\}$ to $\log \zeta$ is already
proved in Lemma \ref{lemma:log}. We need an analogous result for
$\theta(q,\zeta)/\theta(q,cx)$ which is solution of the
$q$-difference equation $y(q\zeta)=cy(\zeta)$, $c\in\C^\ast$. We
give a proof of the needed estimate, although it is a classical
result:

\begin{lemma}\label{lemma:uniformconvergencecar}
Let $c(q)\in\C^\ast$ be a function of $q\in(0,1)$ such that
$\lim_{q\to 1}\frac{\log c}{\ln q}=\gamma$, $\epsilon\in(0,\pi)$
and consider the sector $V_\epsilon:=\{\zeta\in\C^*:
\vert\arg\zeta\vert\le \pi-\epsilon\}$. As $q\to 1^-$, the
following uniform estimate holds uniformly over $V_\epsilon$:
$$
\frac{\theta(q,\zeta)}{\theta(q,c\zeta)}=\zeta^{\gamma}(1+o(1-q))\,.
$$
\end{lemma}

\begin{proof}%[Proof of Lemma \ref{lemma:uniformconvergencecar}]
Let us consider again the modular equation \eqref{equa:Modular} :
\begin{equation*}
\theta(q,\sqrt q\,
)=\sqrt{\frac{2\pi}{\ln(1/q)}}\,e^{-\frac{\log^2x}{2\ln
q}}\theta(q^*,\sqrt{q^*}\,x^*),
\end{equation*}
where $$ x^*=e^{-\frac{2\pi i}{\ln q}\,\log x},\quad
q^*=e^{\frac{4\pi^2}{\ln q}}\,.
$$
We observe that $(\zeta/\sqrt q)^*=-\zeta^*$ and  $(c\zeta/\sqrt
q)^*=-c^*\,\zeta^*$. Therefore we obtain:
$$
\frac{\theta(q,\zeta)}{\theta(q,c\zeta)}=e^{\frac{\log c}{\ln
q}\,(\log(\frac\zeta{\sqrt q})+\frac{\log
c}2)}\,\frac{\theta(q^*,-\sqrt{q^*}\,\zeta^*)}{\theta(q^*,-\sqrt{q^*}\,c^*\,\zeta^*)}\,.
$$
As in the proof of Lemma \ref{lemma:log}, we observe that for any
$\zeta\in V_\epsilon$, $\zeta^*\in V_\epsilon^*$; see
\eqref{equa:Module} for more details. Moreover, when $X$ and $Y$
denote two complex numbers such that $\vert X\vert$, $\vert
Y\vert\in(q^*,1)$, we have the estimate:
$$
\l\vert\frac{\theta(q^*,X)}{\theta(q^*,Y)}\r\vert\le\frac{(-\vert
X\vert;q^*)_\infty\,(-\frac{q^*}{\vert X\vert};q^*)_\infty}{(\vert
Y\vert;q^*)_\infty\,(\frac{q^*}{\vert
Y\vert};q^*)_\infty}\,=\frac{\theta(q^*,\vert
X\vert)}{\theta(q^*,-\vert Y\vert)}\,.
$$
An elementary calculation using \eqref{equa:Module}
allows to conclude, since $q^*\to0^+$ and $c^*\to e^{-2\pi
i\gamma}$.
\end{proof}

Resuming, the function $\phi(q,\xi)$ is a linear combination, with
coefficients that are continuous functions of $q$, of a canonical
basis of solutions at $\infty$: we have proven that both the
canonical basis and the coefficients of the linear combination
admit a uniform limit in a bounded sector containing $d$,
centered at $\infty$ and of arbitrary radius $R>0$. Combined with
Lemma \ref{lemma:uniformconvergenceat0}, this means that
$\phi(q,x)$ converges uniformly to $\phi(1,x)$ in a neighborhood
of the direction $d$, which allows to conclude the proof.

%%%%%%%%%%%%%%%%%%%%%%%%%%%%%%%%%%%%%%%%%%%%%%%%%%%%%%%%%%%%%%%%%%%%
\appendix
%%%%%%%%%%%%%%%%%%%%%%%%%%%%%%%%%%%%%%%%%%%%%%%%%%%%%%%%%%%%%%%%%%%%
\section{Jackson's integral}
\label{sec:jacobsonintegral}
%%%%%%%%%%%%%%%%%%%%%%%%%%%%%%%%%%%%%%%%%%%%%%%%%%%%%%%%%%%%%%%%%%%%

\begin{defn}\label{defn:integral}
We set
$$
F(x)=\int_0^xf(t)\dq t=(1-q)x\sum_{n\geq 0}f(q^n x)q^n\,,
$$
whenever the right hand side converges.
\end{defn}

\begin{rmk}~~
\begin{trivlist}

\item 1.
Notice that if $F(x)$ is well-defined than $\dq F(x)=f(x)$.

\item 2.
If $f(x)$ is continuous on the closed disk $D(0,r^+)$, then $F(x)$
is well defined for any $x\in D(0,r^+)$. In fact there exists
$M>0$ such that $|f(q^n x)q^n|\leq M|q|^n$, which guarantees the
convergence of the infinite sum.

\end{trivlist}
\end{rmk}

\begin{prop}~~
\begin{trivlist}

\item 1.
If $f(x)$ is an analytic function over the disk $D(0,r^-)$, then
$F(x)$ is also analytic over $D(0,r^-)$.

\item 2.
If $F(x)$ is analytic over the disk $D(0,r^-)$ and if $G(x)$ is
another analytic function over $D(0,r^-)$ such that $\dq
G(x)=f(x)$, then
$$
\int_0^x f(t)\dq t=G(x)-G(0)\,.
$$

\end{trivlist}
\end{prop}

\begin{proof}
\begin{trivlist}
\item 1.
It follows from the fact that $F(x)$ is a uniformly convergent
series of analytic function over $D(0,r-\veps^+)$, for any
$r>\veps>0$.

\item 2.
It follows immediately from the remark that the subfield of
constants of the ring of analytic function over $D(0,r^-)$ with
respect to the operator $f(x)\mapsto f(qx)$ is $\C$. In fact this
implies that $F(x)-G(x)\in\C$.
\end{trivlist}
\end{proof}

\begin{defn}\label{defn:integral1}
Let us fix a $q$-orbit $q^{{\mathbb Z}}\a\subset\C$ and suppose
that for any $x\in q^{\mathbb Z}\a$ the integral $F(x)=\int_0^x
f(t)\dq t$ is well-defined. Then we set
$$
\int_{q^Z\a}f(t)\dq t=\lim_{|x|\to\infty\atop x\in q^{\mathbb
Z}\a}\int_0^x f(t)\dq t\,.
$$
\end{defn}

\begin{rmk}\label{rmk:calculintegral}
A straightforward calculation shows that
$$
\int_{q^Z\a}f(t)\dq t=(1-q)\a\sum_{n\in{\mathbb Z}}f(q^n\a)q^n
$$
and in particular that
$$
\int_{q^{\mathbb Z} \a}f(t)\frac{\dq
t}t=(1-q)\sum_{n=-\infty}^{+\infty}f(q^n\a)\,,
$$
whenever the right side converges.
\end{rmk}

%%%%%%%%%%%%%%%%%%%%%%%%%%%%%%%%%%%%%%%%%%%%%%%%%%%%%%%%%%%%%%%%%%
\section{Expansion of $\cE_q(x)$ at $\infty$}
\label{sec:watsonformula}
%%%%%%%%%%%%%%%%%%%%%%%%%%%%%%%%%%%%%%%%%%%%%%%%%%%%%%%%%%%%%%%%%%

The purpose of this section is the proof of Proposition
\ref{prop:qEulerinfinity}. We recall the notation
\begin{equation}\label{equa:an}
a_{n+1}=\sum_{k=0}^n\frac1{q^{k+1}-1},\quad n\ge 0\,,
\end{equation}
$$
A(q)=\sum_{n\ge 0}\frac{q^{n+1}}{q^{n+1}-1}\,,
$$
and
$$
\theta(x)=\sum_{n\in{\mathbb Z}}q^{n(n-1)/2}x^n\,,
$$
and the statement of the proposition:

\begin{prop}
\label{prop:a=q,c=0} For any $x\in\C^*\setminus q^{-{\mathbb N}}$:
\begin{equation}\label{equa:a=q,c=0}
\sum_{n\ge 0}(q;q)_nx^{n+1}=-
\l(-qx\frac{\theta^\prime(-qx)}{\theta(-qx)}+1+A(q)\r)\ \l(\frac
qx;q\r)_\infty+\sum_{n\ge 1}
\frac{a_n}{(q;q)_{n}}q^{n(n+1)/2}\l(-\frac1x\r)^n\,.
\end{equation}
\end{prop}

The proof of the proposition above is based on a Watson's formula
for basic hypergeometric functions, which is the analogue of a
Barnes' formula for Gauss hypergeometric function. Barnes (\cf
\cite[\S 14.51]{WW} and \cite{Barnes}) proved that if $\vert
\arg(-x)\vert <\pi$, $c\notin{\mathbb Z}_{\le 0}$ and
$a-b\notin{\mathbb Z}$, then the analytic continuation of
$_2F_1(a,b;c;x)$ for $\vert x\vert>1$ is given by:
\begin{eqnarray*}
\ds _2F_1(a,b;c;x)&=&\ds\frac{\Gamma(c)\Gamma(b-a)}
    {\Gamma(b)\Gamma(c-a)}(-x)^{-a}{}_2F_1(a,a-c+1,a-b+1;x^{-1})\\
&&\ds +\frac{\Gamma(c)\Gamma(a-b)}{\Gamma(a)\Gamma(c-b)}(-x)^{-b}{}_2F_1(b,b-c+1,b-a+1;x^{-1}).
\end{eqnarray*}
G.N.~Watson (\cf \cite[\S 4.3]{GR}) proved  a formula of the same
kind for Heine series, namely if:
\begin{equation}\label{hyp:abc}
x\notin q^{-{\mathbb N}}\cup \l(\frac{cq}{ab}q^{{\mathbb N}}\r),\quad
c\notin q^{-{\mathbb N}},\quad \frac ab\notin q^{\mathbb Z},\quad
abcx\not=0.
\end{equation}
then
\begin{eqnarray}\label{equa:Watson}
\ds {}_2\phi_1(a,b;c;q,x)&=&
\ds\frac{(b,c/a;q)_\infty}{(c,b/a;q)_\infty}\
\frac{\theta(-ax)}{\theta(-x)}\
{}_2\phi_1\l(a,aq/c;{aq}/b;q,\frac{cq}{abx}\r)\nonumber\\
&&\ds + \frac{(a,c/b;q)_\infty}{(c,a/b;q)_\infty}\
\frac{\theta(-bx)}{\theta(-x)}\
{}_2\phi_1\l(b,{bq}/c;{bq}/a;q,\frac{cq}{abx}\r)\,,
\end{eqnarray}
where $(\alpha_1,...,\alpha_k;q)_n=\prod_{i=1}^k(\alpha_i;q)_n$.
We are going to consider a degeneration of Watson's formula
letting $b\to a$ and $c\to 0$. In this way we obtain an expression
for ${}_2\phi_1(a,a;0;q,x)$ that we can apply to
$$
\cE_q(x)=x{}_2\phi_1\l(q,q;0;q,-\frac x{1-q}\r)\,.
$$

%%%%%%%%%%%%%%%%%%%%%%%%%%%%%%%%%%%%%%%%%%%%%%%%%%%%%%%%%%%%%%%%%%
\subsection{Degenerate cases of the Watson's formula}
\label{subsubsec:watson}
%%%%%%%%%%%%%%%%%%%%%%%%%%%%%%%%%%%%%%%%%%%%%%%%%%%%%%%%%%%%%%%%%%

Let us first consider the case $b\to a q^m$, where $m$ denotes a
non-negative integer. For this purpose, we introduce the following
notation:
$$
\Omega_{m+1}(x)=\sum_{k=0}^m\frac{q^kx}{q^kx-1};
$$
as extension, we set
$$\Omega_0(x)=0,\quad
\Omega(x):=\Omega_\infty(x)=\displaystyle
\lim_{m\to\infty}\Omega_m(x).
$$
Notice that $\Omega_m(x)$ may be identified to a logarithmic
derivative as follows:
$$
\Omega_m(x)=\frac{x\frac{d\ }{dx}(x;q)_m}{(x;q)_m}=x\frac{d\
}{dx}\log(x;q)_m,\quad m\in{\mathbb N}\cup\{\infty\}.
$$
Since $(x;q)_{m+n}=(x;q)_m\,(xq^m;q)_n$ for any $m$, $n\in{\mathbb
N}$, it follows:
\begin{equation}\label{equa:Omega}
 \Omega_{m+n}(x)=\Omega_m(x)+\Omega_n(q^mx),\quad \Omega(x)=\Omega_m(x)+\Omega(q^mx).
\end{equation}
Let
$$
\Gamma_q(x)=\frac{(q;q)_\infty}{(q^x;q)_\infty}(1-q)^{1-x}\,,
\hbox{ defined for any $x\notin(-\N)$,}
$$
be the Jackson's Gamma function (\cf \cite[\S10.1]{GR}).
It's useful to remark that, if
$\ds\Psi_q(x)=\frac{d\ }{dx}\log \Gamma_q(x)$, then
\begin{equation}\label{equa:PsiOmega}\Psi_q(x)=-\ln q\,\Omega(q^x)-\ln(1-q),\quad
\Psi_q(x+m)=\Psi_q(x)+\ln q\,\Omega_m(q^x)\,,
\end{equation}
for any non-negative integer $m$.
\par
Let us define
$$
\ell(x):=\ell_q(x)=-x\frac{\theta'(-x)}{\theta(-x)}\,,
$$
where $\theta'(-x)$ denotes the derivative of $\theta$ w.r.t. at
the variable $x$. From the  Jacobi's triple formula
$\theta(x)=(q,-x,-q/x;q)_\infty$, one deduces the following
relation:
\begin{equation}\label{equa:qlogOmega}
\ell(x)=-\Omega(x)+\Omega(\frac{q}{x}).
\end{equation}
Since $\Omega_1(x)+\Omega_1(\frac1x)=1$, putting $m=1$ in
$\eqref{equa:Omega}$ allows to obtain the following relation:
$$
\ell(qx)-\ell(x)=1,
$$
which means that, for any given non-zero complex number $\la$, the
function $x\mapsto\ell(\la x)$ is a $q$-logarithm. From
\eqref{equa:PsiOmega}, one gets the following link between
$\ell_q$ and $\Psi_q$:
\begin{equation}\label{equa:qlogPsi}
\ell(q^x)=\frac{1}{\ln q}\, \l(\Psi_q(x)-\Psi_q(1-x)\r)\,.
\end{equation}

\begin{prop}\label{prop:b=aq^m}
Let $m$ be a non-negative integer and let $a$, $c$ be non-zero
complex numbers. Suppose that $a\notin q^{-{\mathbb N}}$, $c\notin
q^{-{\mathbb N}}$ and $c/a\notin q^{-{\mathbb N}}$. Then, the
following formula holds:
{\setlength\arraycolsep{2pt}
\begin{eqnarray}\label{equa:b=aq^m}
{}_2\phi_1(a,aq^m;c;q,x)&=&\frac{(aq^m,c/a;q)_\infty}{(c,q^m;q)_\infty}\, \frac{\theta(-ax)}{\theta(-x)}\,P_m(a,c,x)+\frac{(a,{cq^{-m}}/{a};q)_\infty}{(c,q;q)_\infty\,(q^{-m};q)_m}\,\frac{\theta(-aq^mx)}{\theta(-x)}\, \nonumber \\
&&\times\l\{\l[C_m(a,c)+\ell(aq^mx)\r]\Phi_m(a,c,x)+\sum_{n=1}^\infty
C_{m,n}(a,c)\,\phi_{m,n}(a,c,x)\r\} \,,
\end{eqnarray}}
where
$$
P_m(a,c,x)=\frac{(q;q)_{m-1}}{(a;q)_m}\,
\sum_{n=0}^{m-1}\frac{(a,{aq}/{c};q)_n}{(q,q^{1-m};q)_n}\,
\l(\frac{cq^{1-m}}{a^2x}\r)^n,
$$
$$
C_m(a,c)=\Omega(q)+\Omega(q^{1+m})-\Omega(cq^{-m}/a)-\Omega(aq^m)+1,
$$
$$
\Phi_m(a,c,x)={}_2\phi_1(aq^m,aq^{1+m}/c;q^{1+m};q,\frac{cq^{1-m}}{a^2x})=\sum_{n\ge
0}\phi_{m,n}(a,c,x) ,
$$
$$
\phi_{m,n}(a,c,x)=\frac{(aq^m,aq^{1+m}/c;q)_n}{(q,q^{1+m};q)_n}\,\l(\frac{cq^{1-m}}{a^2x}\r)^n
$$
and
$$
C_{m,n}(a,c)=\Omega_n(aq^m)+\Omega_n(aq^{1+m}/c)-\Omega_n(q^{1+m})-\Omega_n(q).
$$
When $m=0$, $P_m(a,c,x)=0$.
\end{prop}

\begin{rmk}
Equations \eqref{equa:b=aq^m} is a $q$-analog of \cite[p. 109, (7)]{erd}.
\end{rmk}

\begin{proof}
Letting $b=aq^m\epsilon$ in \eqref{equa:Watson} gives raise to the
following formula:
{\setlength\arraycolsep{2pt}
\begin{eqnarray}
\label{equa:Watson1}
{}_2\phi_1(&a&,aq^m\epsilon;c;q,x)=\frac{(aq^m\epsilon,c/a;q)_\infty}{(c,q^m\epsilon;q)_\infty}\,
\frac{\theta(-ax)}{\theta(-x)}\,
{}_2\phi_1\l(a,aq/c;q^{1-m}/\epsilon;q,\frac{cq^{1-m}}{a^2\epsilon x}\r)\nonumber\\
&&+
\frac{(a,cq^{-m}/(a\epsilon);q)_\infty}{(c,q^{-m}/\epsilon;q)_\infty}\
\frac{\theta(-aq^m\epsilon x)}{\theta(-x)}\
{}_2\phi_1\l(aq^m\epsilon,{aq^{1+m}\epsilon}/c;q^{1+m}\epsilon;q,\frac{cq^{1-m}}{a^2\epsilon
x}\r)\,,
\end{eqnarray}}
Suppose $m\ge 1$. Since
$$
(X;q)_{n+m}=(X;q)_m\,(Xq^m;q)_n,\quad
(X;q)_m=(-X)^m\,(q^{1-m}/X;q)_m\,q^{m(m-1)/2}
$$
and
$$
\theta(q^mX)=X^{-m}\,q^{-m(m-1)/2}\,\theta(X),
$$
we obtain
$$
\begin{array}{rcl}
\lefteqn{\hskip -15pt{}_2\phi_1\l(a,aq/c;q^{1-m}/\epsilon;q,\frac{cq^{1-m}}{a^2\epsilon x}\r)
=\sum_{n=0}^{m-1}\frac{(a,aq/c;q)_n}{(q,q^{1-m}/\epsilon;q)_n}
\,\l(\frac{cq^{1-m}}{a^2\epsilon x}\r)^n}\\
&+&\ds
\frac{(a,cq^{-m}/a;q)_m}{(q,\epsilon;q)_m}\,\l(\frac{q}{a\epsilon
x}\r)^m\,\sum_{n=0}^\infty\frac{(aq^m,aq^{1+m}/c;q)_n}{(q^{1+m},q/\epsilon;q)_n}\,\l(\frac{cq^{1-m}}{a^2\epsilon
x}\r)^n
\end{array}
$$
and
$$
\frac{(a,cq^{-m}/(a\epsilon);q)_\infty}{(c,q^{-m}/\epsilon;q)_\infty}\
\frac{\theta(-aq^m\epsilon
x)}{\theta(-x)}=-\epsilon\,\frac{(a,cq^{-m}/(a\epsilon);q)_\infty}{(c,q/\epsilon;q)_\infty\,
(\epsilon;q)_{m+1}}\,\l(\frac q{ax}\r)^m\,\frac{\theta(-a\epsilon x)}{\theta(-x)}\,.
$$
Hence, we can re-write \eqref{equa:Watson1} as follows:
\begin{equation}\label{equa:Watson2}
{}_2\phi_1(a,aq^m\epsilon;c;q,x)=A(\epsilon)+\frac{B_1(\epsilon)C_1(\epsilon)-B_2(\epsilon)C_2(\epsilon)}{\epsilon-1
}\,,
\end{equation}
where
$$A(\epsilon)=
\frac{(aq^m\epsilon,c/a;q)_\infty}{(c,q^m\epsilon;q)_\infty}\,
\frac{\theta(-ax)}{\theta(-x)}\,
\sum_{n=0}^{m-1}\frac{(a,aq/c;q)_n}{(q,q^{1-m}/\epsilon;q)_n}\,\l(\frac{cq^{1-m}}{a^2\epsilon
x}\r)^n\,,
$$
$$
B_1(\epsilon)=\epsilon\,\frac{(a,cq^{-m}/(a\epsilon);q)_\infty}{(c,q/\epsilon;q)_\infty\,(\epsilon
q;q)_m }\,\l(\frac q{ax}\r)^m\,\frac{\theta(-a\epsilon
x)}{\theta(-x)}\,,
$$
$$
C_1(\epsilon)={}_2\phi_1
\l(aq^m\epsilon,{aq^{1+m}\epsilon}/c;q^{1+m}
\epsilon;q,\frac{cq^{1-m}}{a^2\epsilon x}\r)\,,
$$
$$
B_2(\epsilon)=\frac{(a,cq^{-m}/a;q)_\infty}{(c,q\epsilon;q)_\infty\,(q;q)_m}\,\frac{(aq^{m}\epsilon;q)_\infty}{(aq^m;q)_\infty}\,\l(\frac{q}{a\epsilon
x}\r)^m\, \frac{\theta(-ax)}{\theta(-x)}\,,
$$
$$
C_2(\epsilon)=\sum_{n=0}^\infty\frac{(aq^m,aq^{1+m}/c;q)_n}{(q^{1+m},q/\epsilon;q)_n}\,\l(\frac{cq^{1-m}}{a^2\epsilon
x}\r)^n\,.
$$
Since $ B_1(1)=B_2(1)$ and $C_1(1)=C_2(1)$, putting $\epsilon\to
1$ in \eqref{equa:Watson2} allows us to get the following
relation:
$$
{}_2\phi_1(a,aq^m;c;q,x)=A(1)+\l[B_1'(1)-B_2'(1)\r]C+B\l[C_1'(1)-C_2'(1)\r]\,,
$$
with $C=C_1(1)$, $B=B_1(1)$. By direct computation, it yields:
$$
B_1'(1)=\l[1-\Omega(cq^{-m}/a)+\Omega(q)-\Omega_m(q)+\ell(ax)\r]\,B,\quad
B_2'(1)=\l[\Omega(aq^m)-\Omega(q)-m\r]\,B,
$$
$$
C_1'(1)=\sum_{n=1}^\infty\frac{(aq^m,aq^{1+m}/c;q)_n}{(q^{1+m},q;q)_n}\,
\l[\Omega_n(aq^m)+\Omega_n(aq^{1+m}/c)-\Omega_n(q^{1+m})-n\r]\,\l(\frac{cq^{1-m}}{a^2
x}\r)^n
$$
and
$$
C_2'(1)=\sum_{n=1}^\infty\frac{(aq^m,aq^{1+m}/c;q)_n}{(q^{1+m},q;q)_n}\,
\l[\Omega_n(q)-n\r]\,\l(\frac{cq^{1-m}}{a^2 x}\r)^n\,.
$$
Notice also that
$$
B=\frac{(a,cq^{-m}/a;q)_\infty}{(c,q;q)_\infty\,( q;q)_m
}\,\l(\frac q{ax}\r)^m\,\frac{\theta(-a
x)}{\theta(-x)}=\frac{(a,cq^{-m}/a;q)_\infty}{(c,q;q)_\infty\,(
q^{-m};q)_m }\,\frac{\theta(-aq^m x)}{\theta(-x)}\,,
$$
which ends the proof when $m\ge 1$.

If $m=0$, $A(\epsilon)$ doesn't exist, so $P_0(a,c,x)=0$.
\end{proof}

Consider now two cases for which the hypothesis of Proposition
\ref{prop:b=aq^m} are not all fulfilled: $c/a\in q^{-{\mathbb N}}$
or $c=0$. Let $k\in{\mathbb N}$ and $c/a=q^{-k}\epsilon$, where
$\epsilon\to 1$. Since
$$\lim_{\epsilon\to 1}(q^{-k}\epsilon;q)_\infty\,\Omega(q^{-m-k}\epsilon)=-(q^{-k};q)_k\,(q;q)_\infty,
$$
it follows:

\begin{cor}
Let $k$, $m\in{\mathbb N}$ and $a\in\C^*$. The following equality
holds:
\begin{equation}\label{equa:b=aq^ma=cq^k}
{}_2\phi_1(a,aq^m;aq^{-k};q,x)=
\frac{(q^{-k-m};q)_k}{(aq^{-k};q)_k}\,\frac{\theta(-ax)}{\theta(-x)}\,
{}_2\phi_1\l(aq^m,q^{1+m+k};q^{1+m};q,\frac{q^{1-k-m}}{ax}\r)\,.
\end{equation}
\end{cor}

Taking now $c\to 0$ in \eqref{equa:b=aq^m} gives the following
formula:
{\setlength\arraycolsep{2pt}\begin{eqnarray}\label{equa:b=aq^mc=0}
{}_2\phi_1(a,aq^m;0;q,x)&=&\frac{(aq^m;q)_\infty}{(q^m;q)_\infty}\, \frac{\theta(-ax)}{\theta(-x)}\,P_m(a,x)+\frac{(a;q)_\infty}{(q;q)_\infty\,(q^{-m};q)_m}\,\frac{\theta(-aq^mx)}{\theta(-x)}\, \nonumber \\
&&\times\l\{\l[C_m(a)+\ell(aq^mx)\r]\Phi_m(a,x)+\sum_{n=1}^\infty
C_{m,n}(a)\,\phi_{m,n}(a,x)\r\} \,,
\end{eqnarray}}
where
$$
P_m(a,x)=\frac{(q;q)_{m-1}}{(a;q)_m}\,\sum_{n=0}^{m-1}\frac{(a;q)_n\,q^{n(n-1)/2}}{(q,q^{1-m};q)_n}\,\l(-\frac{q^{2-m}}{ax}\r)^n,
$$
$$
C_m(a)=\Omega(q)+\Omega(q^{1+m})-\Omega(aq^m)+1,
$$
$$
\Phi_m(a,x)={}_1\phi_1(aq^m;q^{1+m};q,\frac{q^{2}}{ax})=\sum_{n\ge
0}\phi_{m,n}(a,x) ,
$$
$$
\phi_{m,n}(a,x)=\frac{(aq^m;q)_n\,q^{n(n-1)/2}}{(q,q^{1+m};q)_n}\,\l(-\frac{q^{2}}{ax}\r)^n
$$
and
$$
C_{m,n}(a)=\Omega_n(aq^m)+n-\Omega_n(q^{1+m})-\Omega_n(q).
$$
Once again, when $m=0$ and $P_m(a,x)=0$, we have:
\begin{equation}\label{equa:b=ac=0}
{}_2\phi_1(a,a;0;q,x)=\frac{(a;q)_\infty}{(q;q)_\infty}\,\frac{\theta(-ax)}{\theta(-x)}\,
\l\{\l[C_0(a)+\ell(ax)\r]\Phi_0(a,x)+\sum_{n=1}^\infty
C_{0,n}(a)\,\phi_{0,n}(a,x)\r\} \,,
\end{equation}
where
$$
C_0(a)=2\Omega(q)-\Omega(a)+1,\quad
\Phi_0(a,x)={}_1\phi_1(a;q;q,\frac{q^{2}}{ax})=
\sum_{n\ge 0}\phi_{0,n}(a,x) ,
$$
$$
\phi_{0,n}(a,x)=\frac{(a;q)_n\,q^{n(n-1)/2}}{(q,q;q)_n}\,\l(-\frac{q^{2}}{ax}\r)^n
\,,\quad C_{0,n}(a)=\Omega_n(a)+n-2\Omega_n(q).
$$

%%%%%%%%%%%%%%%%%%%%%%%%%%%%%%%%%%%%%%%%%%%%%%%%%%%%%%%%%%%%%%%%%%
\subsection{Proof of Proposition
\ref{prop:qEulerinfinity}}\label{subsec:ProofProp}
%%%%%%%%%%%%%%%%%%%%%%%%%%%%%%%%%%%%%%%%%%%%%%%%%%%%%%%%%%%%%%%%%%%%

The equality ${}_1\phi_1(q;q;q,X)=(X;q)_\infty$,
plus \eqref{equa:b=ac=0}, where we have set $a=q$, implies the following
formula:
$$
\sum_{n\ge 0}(q;q)_nx^n=-
\frac{1}{x}\l\{\l[\Omega(q)+1+\ell(qx)\r]\,(q/x;q)_\infty+\sum_{n\ge
1}\l[n-\Omega_n(q)\r]\frac{q^{n(n-1)/2}}{(q;q)_n}\,\l(-\frac
qx\r)^n\r\}.
$$
Thus, one may obtain Proposition \ref{prop:qEulerinfinity} by
taking into account the following relations:
$$
\ell(qx)=-qx\,\frac{\theta'(-qx)}{\theta(-qx)},\quad
a_{n}=\Omega_{n}(q)-n, \quad A_q=\Omega(q).
$$

%\pagebreak

%%%%%%%%%%%%%%%%%%%%%%%%%%%%%%%%%%%%%%%%%%%%%%%%%%%%%%%%%%%%%%%%%%%%%%%%%%%%%%%%%%%%%%%%
%%%%%%%%%%%%%%%%%%%%%%%%%%%%%%%%%%%%%%%%%%%%%%%%%%%%%%%%%%%%%%%%%%%%%%%%%%%%%%%%%%%%%%%%
%%%%%%%%%%%%%%%%%%%%%%%%%%%%%%%%%%%%%%%%%%%%%%%%%%%%%%%%%%%%%%%%%%%%%%%%%%%%%%%%%%%%%%%%
%%%%%%%%%%%%%%%%%%%%%%%%%%%%%%%%%%%%%%%%%%%%%%%%%%%%%%%%%%%%%%%%%%%%%%%%%%%%%%%%%%%%%%%%
%%%%%%%%%%%%%%%%%%%%%%%%%%%%%%%%%%%%%%%%%%%%%%%%%%%%%%%%%%%%%%%%%%%%%%%%%%%%%%%%%%%%%%%%
%%%%%%%%%%%%%%%%%%%%%%%%%%%%%%%%%%%%%%%%%%%%%%%%%%%%%%%%%%%%%%%%%%%%%%%%%%%%%%%%%%%%%%%%
%%%%%%%%%%%%%%%%%%%%%%%%%%%%%%%%%%%%%%%%%%%%%%%%%%%%%%%%%%%%%%%%%%%%%%%%%%%%%%%%%%%%%%%%
%%%%%%%%%%%%%%%%%%%%%%%%%%%%%%%%%%%%%%%%%%%%%%%%%%%%%%%%%%%%%%%%%%%%%%%%%%%%%%%%%%%%%%%%

\section*{Part II. Summation of divergent $q$-series and confluence: the case $q>1$}
\addcontentsline{toc}{section}{Part II. Summation of divergent
$q$-series and confluence: the case $q>1$}
%%%%%%%%%%%%%%%%%%%%%%%%%%%%%%%%%%%%%%%%%%%%%%%%%%%%%%%%%%%%%%%%%%%%%%%%%%%%%%%%%%%%%%%%
%%%%%%%%%%%%%%%%%%%%%%%%%%%%%%%%%%%%%%%%%%%%%%%%%%%%%%%%%%%%%%%%%%%%%%%%%%%%%%%%%%%%%%%%
%%%%%%%%%%%%%%%%%%%%%%%%%%%%%%%%%%%%%%%%%%%%%%%%%%%%%%%%%%%%%%%%%%%%%%%%%%%%%%%%%%%%%%%%
%%%%%%%%%%%%%%%%%%%%%%%%%%%%%%%%%%%%%%%%%%%%%%%%%%%%%%%%%%%%%%%%%%%%%%%%%%%%%%%%%%%%%%%%

\makeatletter
\renewcommand{\thesection}{\@arabic\c@section}
\makeatother

{\bfseries Important.} From now on, we fix
$q\in(1,+\infty)\subset{\mathbb R}$, so that $p=q^{-1}\in(0,1)$.

\bigskip
In this second part we introduce four types of $q$-summation (\cf
Definition \ref{defn:summations} below): our purpose is studying
the relations among them. First of all, we investigate the
different sums of the $q$-Euler series $\sum_{n\geq
0}(-1)^n[n]_q^!x^{n+1}$ and their properties (\cf \S
\ref{sec:qEulersums} below). Then we prove a general result for
generic $q$-Gevrey series (\cf Theorem \ref{thm:sums}), based on
the study of the Tschakaloff series
$$
T_q(x)=\sum_{n\geq 0}q^{n(n-1)/2}x^{n+1}\,,
$$
and of a convenient $q$-convolution product.

\medskip\noindent
{\bf Notation.} We set:
\begin{equation}\label{eq:Jacobi}
\theta_p(x):=\theta(p,x)=\sum_{n\in{\mathbb Z}}p^{n(n-1)/2}x^n
=(p;p)_\infty(-x;p)_\infty(-p/x;p)_\infty\,
\end{equation}
and
\begin{equation}\label{eq:q-exp}
e_q(x):=(-(1-p)x;p)_\infty=\sum_{n\geq 0}\frac{x^n}{[n]_q^!}\,.
\end{equation}
Remark that $e_q(x)=e_p(-x)^{-1}$ and that
$\theta_p(x)=x\theta_p(px)$.

%%%%%%%%%%%%%%%%%%%%%%%%%%%%%%%%%%%%%%%%%%%%%%%%%%%%%%%%%%%%%%%%%%%%%%%%%%%%%%%%%%%%%%%%
%%%%%%%%%%%%%%%%%%%%%%%%%%%%%%%%%%%%%%%%%%%%%%%%%%%%%%%%%%%%%%%%%%%%%%%%%%%%%%%%%%%%%%%%
\section{The divergent $q$-Euler series}
\label{sec:qEulersums}
%%%%%%%%%%%%%%%%%%%%%%%%%%%%%%%%%%%%%%%%%%%%%%%%%%%%%%%%%%%%%%%%%%%%%%%%%%%%%%%%%%%%%%%%
%%%%%%%%%%%%%%%%%%%%%%%%%%%%%%%%%%%%%%%%%%%%%%%%%%%%%%%%%%%%%%%%%%%%%%%%%%%%%%%%%%%%%%%%

Since $q>1$, the $q$-Euler series
$$
\hat\cE_q(x)=\sum_{n\geq 0}(-1)^n[n]_q^!x^{n+1}
$$
is obviously divergent for any $x\in\C^\ast$, as the Euler series
$\sum_{n\ge 0}(-1)^n n!x^{n+1}$. The corresponding $q$-difference
equation is
$$
x^2d_qy+y=x\,.
$$

%%%%%%%%%%%%%%%%%%%%%%%%%%%%%%%%%%%%%%%%%%%%%%%%%%%%%%%%%%%%%%%%%%%%%%%%%%%%%%%%%%%%%%%%
\subsection{Definition of different sums of the $q$-Euler series}
\label{subs:sumsEuler}
%%%%%%%%%%%%%%%%%%%%%%%%%%%%%%%%%%%%%%%%%%%%%%%%%%%%%%%%%%%%%%%%%%%%%%%%%%%%%%%%%%%%%%%%

Let us consider the \emph{$q$-Borel transforms of $\hat\cE_q(x)$}
(for the general definition,\cf \S\ref{subsec:borel}):
$$
\psi(\xi):=\frac 1{1+\xi} \hbox{ and } \phi(\xi):=\cE_p(\xi)\,.
$$
In the following, we will identify $\cE_p(\xi)$ to its analytic
continuation on $\C\setminus((p-1)q^{\mathbb N})$. For any
$d\in(-\pi,\pi)$ and $\la\notin -p^{\mathbb Z}$, we set:
\begin{equation}\label{equa:cEd}
\cE_q^d(x)=\frac{q-1}{\ln q}
    \int_0^{e^{id}\infty}\frac{\psi(\xi)}{e_q(q\frac \xi x)}d\xi,\quad \arg x\in(d-\pi,d+\pi);
\end{equation}
\begin{equation}\label{equa:cElambda}
\cE_q^{[\lambda]}(x)=\frac q{1-p}
    \int_{\la p^{\mathbb Z}}\frac{\psi\l(\frac \xi{1-p}\r)}
    {e_q\l(q\frac \xi{(1-p)x}\r)}\,{d_p\xi},\quad x\notin (p-1)\lambda q^{\mathbb Z};
\end{equation}
\begin{equation}\label{equa:Ed}
 E_q^d(x)=\frac{q}{\ln q}
    \int_0^{e^{id}\infty}\frac{\phi(\xi)}{\theta_p(q\frac \xi x)}d\xi,\quad \arg x\in(d-\pi,d+\pi);
\end{equation}
\begin{equation}\label{equa:Elambda}
E_q^{[\lambda]}(x)=\frac q{1-p}
    \int_{\la p^{\mathbb Z}}\frac{\phi(\xi)}{\theta_p(q\frac \xi x)}\,{d_p\xi},\quad x\notin (p-1)\lambda q^{\mathbb Z};
\end{equation}

\begin{prop}\label{prop:cEEd}\begin{trivlist}
 \item
(1) The functions $\cE_q^d$ and $E_q^d$ can be analytically
continued on the sector $\{|\arg x|<2\pi\}$ of the Riemann surface
of the logarithm.
\item
(2) The functions $\cE_q^{[\la]}$ and $E_q^{[\la]}$ are analytic
on $\C^*\setminus (-(1-p)\la q^{\mathbb Z})$, the point $-(1-p)\la
q^n$ being a simple pole for any integer $n\in{\mathbb Z}$.
\end{trivlist}
\end{prop}

\begin{proof}
The functions $\cE_q^d$ and $E_q^d$ are \emph{a priori} defined
for $|\arg x-d|<\pi$ and $d$ varies in $(-\pi,\pi)$. The second
assertion is straightforward.
\end{proof}

We will denote by $\cE_q$ and $E_q$ the analytic continuation of
$\cE_q^d$ and $E_q^d$, respectively, on the open sector
$V(-2\pi,2\pi):=\{x\in\tilde\C^*: \vert\arg x\vert<2\pi\}$ on the
Riemann surface of the logarithm. We recall the following result:

\begin{prop}
[{\cite[Thm. 2.1]{RZ} and \cite[Thm. 1.3.2]{zhanggroningen}}] The
function $\cE_q(x)$ (resp. $\cE_q^{[\la]}$) admits $\hat\cE_q(x)$
as $q$-Gevrey asymptotic expansion at $x=0$ in the sector $\{\arg
x<3\pi/2\}$. In particular they are solution of $x^2\dq y+y=x$.
\end{prop}

The following theorem is about the comparison between the four
summations of $\hat\cE_q$ we have just introduced:

\begin{thm}\label{thm:identity}
$\cE_q(x)=E_q(x)$ and $\cE^{[\la]}_q(x)=E^{[\la]}_q(x)$.
\end{thm}

First, we need to prove the following two lemmas.

\begin{lemma}\label{lemma:stokes}
For any $x\in\C$ such that $\arg x\in(-2\pi,0)$ we have:
\begin{equation}\label{equa:EdStokes}
\cE_q(xe^{2\pi i})-\cE_q(x) =E_q(x e^{2\pi i})-E_q(x) =-2\pi
i\frac{q-1}{\ln q}\frac{1}{e_q(-\frac{q}{x})}\,.
\end{equation}
In particular: $\cE_q(xe^{2\pi i})-E_q(x e^{2\pi
i})=\cE_q(x)-E_q(x)$.
\end{lemma}

\begin{proof}
A variable change in the integral defining $\cE_q^d$ (resp.
$E_q^d$) shows that $\cE_q^d(xe^{2\pi i})=\cE_q^{d-2\pi}(x)$
(resp. $E^d_q(x e^{2\pi i})=E^{d-2\pi}_q(x)$). We are reduced to
calculate $\cE_q^{d-2\pi}(x)-\cE_q^d(x)$ and hence to calculate a
residue at $\xi=-1$. In an analogous way, using formula
\eqref{equa:Heine}, we obtain
$$
\begin{array}{rcl}
E^d_q(x e^{2\pi i})-E^d_q(x) &=&\ds-\frac{q2\pi i}{\ln
q}\sum_{n\ge 0}
    {\rm Res}_{\xi=(p-1)q^n}\l(\frac{\cE_p(\xi)}{\theta_p\l(q\frac{\xi}{x}\r)}\r)\\ \\
&=&\ds -\frac{q(p-1)2\pi i}{\ln q}\sum_{n\ge
0}\frac{q^n(p^{n+1};p)_\infty}
    {\theta_p\l(q\frac{(p-1)q^n}{x}\r)}\\ \\
&=&\ds\frac{q(p-1)2\pi i}{\ln
q}\frac{(p;p)_\infty}{\theta_p\l(\frac{1-q}{x}\r)}
    \sum_{n\ge 0}\frac{p^{n(n-1)/2}}{(p;p)_n}\\ \\
&=&\ds -\frac{q(p-1)2\pi i}{\ln q}
    \frac{\l(p,-\frac x{q-1};p\r)_\infty}{\theta_p\l(\frac{1-q}{x}\r)}\,.
\end{array}
$$
The Jacobi triple product formula for $\theta_p$ immediately
allows to conclude.\end{proof}

\begin{lemma}\label{lemma:homogene}
Let us consider the homogenous $q$-difference equation
\begin{equation}\label{equa:homogene}
x^2d_qy=y.
\end{equation}
Let $y_0$ be a meromorphic solution of \eqref{equa:homogene} on
the domain $\Omega=\{0<\vert x\vert<R\}$. Suppose that one of the
following hypotheses is verified:
\begin{itemize}
\item
the function $y_0$ is analytic on $\Omega$;

\item
there exists $\mu\in\C^*$ such that $\mu\notin(1-p)p^{\mathbb N}$
and such that the function $y_0$ has only simple poles contained
in the set $\mu p^{\mathbb N}$;
\end{itemize}
then $y_0$ is identically zero.
\end{lemma}

\begin{proof}
Notice that $1/e_p(q/x)$ is a uniform solution to
\eqref{equa:homogene}. Hence, there exists a $q$-invariant
function $K(x)$ such that $y_0(x)=K(x)/e_p(q/x)$. Identifying
$K(x)$ to an elliptic function, one ends the proof noticing that
$(1-p)p^{\mathbb Z}$ is the only spiral of poles of $e_p(q/x)$.
\end{proof}

\begin{proof}
[Proof of Theorem \ref{thm:identity}] Lemma \ref{lemma:stokes}
implies that $h^d(x):=E^d_q(x)-\cE^d_q(x)$ is an analytic solution
of \eqref{equa:homogene} on $\C^*$. We deduce from Lemma
\ref{lemma:homogene} that  $h^d\equiv0$.
\par
The difference $E^{[\la]}_q(x)-\cE^{[\la]}_q(x)$ has only simple
poles on $-\lambda(1-p)q^{\mathbb Z}$. Since $\lambda\notin
-p^{\mathbb Z}$ we conclude applying Lemma \ref{lemma:homogene}.
\end{proof}

%%%%%%%%%%%%%%%%%%%%%%%%%%%%%%%%%%%%%%%%%%%%%%%%%%%%%%%%%%%%%%%%%%%%%%%%%%%%%%%%%%%%%%%%
\subsection{$q$-integral and continuous integral}
%%%%%%%%%%%%%%%%%%%%%%%%%%%%%%%%%%%%%%%%%%%%%%%%%%%%%%%%%%%%%%%%%%%%%%%%%%%%%%%%%%%%%%%%

%\begin{thm}\label{thm:integrals}
%$$
%\cE_q(x)-\cE^{[\la]}_q(x)
%=-2\pi i\frac{q-1}{\ln q}x^\ast
%\frac{\theta_{p^\ast}^\p(x^\ast)}{\theta_{p^\ast}^\p(x^\ast)}\, e_q(-q/x)\,,
%\,\forall \arg x\in(-\pi,\pi)\,,
%$$
%where $x^*=e^{2\pi i\log_qx}$ and $q^*=e^{4\pi^2/\ln q}=(p^\ast)^{-1}$.
%\end{thm}
%
%\begin{proof}
%The function $\cE_q^{d}-\cE_q^{[\la]}$ is solution of
%the homogeneous equation $x^2d_qy=-y$, therefore
%$z(x)=\l(\cE_q^{d}-\cE_q^{[\la]}\r)e_q(-q/x)$ verifies the functional equations:
%$$
%z(xe^{2\pi i})-z(x)=-2\pi i\frac{q-1}{\ln q},\qquad z(qx)=z(x)\,.
%$$
%Consider the variable change $x^*=e^{2\pi i\log_qx}$ and set
%$q^*=e^{4\pi^2/\ln q}=(p^\ast)^{-1}$. A direct calculation shows that
%$z(q^\ast x^\ast)-z(x^\ast)=-2\pi i\frac{q-1}{\ln q}$
%and hence that
%$$
%z(x^\ast)=-2\pi i\frac{q-1}{\ln q}x^\ast
%\frac{\theta_{p^\ast}^\p(x^\ast)}{\theta_{p^\ast}(x^\ast)}
%+(\hbox{$q^\ast$-invariant function})\,.
%$$
%Once again the control of the poles of $z(x)$ allows to conclude that
%$z(x^\ast)=-2\pi i\frac{q-1}{\ln q}x^\ast
%\frac{\theta_{p^\ast}^\p(x^\ast)}{\theta_{p^\ast}(x^\ast)}$.
%\end{proof}

Although both $\cE_q(x)$ and $\cE_q^{[\la]}$ are solutions of the
$q$-difference equation $x^2d_qy+y=x$, they have a deeply
different nature. In fact, while $\cE_q(x)$ is meromorphic on the
whole Riemann surface $\wtilde\C^\ast$ of the logarithme, the
function $\cE_q^{[\la]}$ is a uniform function: more precisely, it
is analytic on $\C^\ast$ minus a spiral of simple poles.
\par
Let us consider the projection:
$$
\begin{array}{rccc}
\pi: & ]0,+\infty[\times{\mathbb R}\cong~\wtilde\C^\ast  & \longrightarrow & \C^\ast\\ \\
& (r,\a) & \longmapsto & re^{2i\pi\a}
\end{array}\,.
$$
Of course, we can identify $\cE_q^{[\la]}$ to its pull back
\emph{via} $\pi$ on $\wtilde\C^\ast$, \ie to a meromorphic
function on $\wtilde\C^\ast$, and study the solution
$\cE_q(x)-\cE_q^{[\la]}$ of $x^2\dq y+y=0$. We have the following
result (we identify all meromorphic function on $\C^\ast$ to their
pull back on $\wtilde\C^\ast$):

\begin{prop}\label{prop:disc-cont-euler}
Let $\la\in\C^*\setminus(-q^{\mathbb Z})$. For any
$x\in\tilde\Omega_\la:=\pi^{-1}\l(\C^*\setminus[-\la(1-p);q]\r)$,
we have:
$$
\cE_q(x)-\cE_q^{[\la]}(x)=-2\pi i\frac{q-1}{\ln q}
\l(L_{(-\la(1-p))^*,q^*}(x^*)-C_{\la,q}\r)e_q(-q/x)^{-1},
$$
where the following notations are employed:
$$
\begin{array}{ll}
\ds x^*=e^{-2\pi i\frac{\log x}{\ln q}}\,,
&\ds q^*=e^{4\pi^2/\ln q}\,,\\ \\
\ds L_{a,q}(x)=-\frac xa\frac{\theta_p'(-\frac
xa)}{\theta_p(-\frac xa)}=\ell_p\l(\frac xa\r)\,, &\ds
C_{\la,q}=L_{(-\la(1-p))^*,q^*}\l((1-p)^*\r)\,.
\end{array}
$$
\end{prop}

Since $L_{a,q}(\frac x{a'})=L_{aa',q}(x)=\ell_p(\frac x{aa'})$ and
$(xx')^*=x^*(x')^*$, the theorem above can be rephrased in the
following statement:

\begin{cor}\label{cor:disc-cont-euler} Let
$p^*=1/q^*=e^{-4\pi^2/\ln q}$. Then
$$
\cE_q(x)-\cE_q^{[\la]}(x)=-2\pi i\frac{q-1}{\ln q}
\l[\ell_{p^*}\l(\Big(-\frac x{\la(1-p)}\Big)^*\r)-
\ell_{p^*}\l(\Big(-\frac1\la\Big)^*\r)\r]
e_q(-q/x)^{-1}\,.
$$
\end{cor}

The proof  of Proposition \ref{prop:disc-cont-euler}
is based on the following two lemmas:

\begin{lemma}\label{lemma:U}
Let $\la\in\C^*\setminus(-p^{\mathbb Z})$ and
$\tilde\Omega_\la:=\pi^{-1}(\C^*\setminus[-\la(1-p);q])\subset\tilde\C^*$.
For any $x\in\tilde\Omega_\la$ we set:
$$
U_\la(x)=(\cE_q(x)-\cE^{[\la]}(x))e_q(-q/x)\,.
$$
Then
$$
U_\la(xe^{2\pi i})-U_\la(x)=-2\pi i\frac{q-1}{\ln q}
$$
and
$$
U_\la(qx)=U_\la(x)\,.
$$
\end{lemma}

\begin{proof}
The proof follows from Lemma \ref{lemma:stokes}, taking into
account the functional equation of $e_q(-q/x)$:
$$
e_q(-q/x)=\l(1-\frac{q-1}x\r)e_q(-1/x)
$$
\end{proof}

\begin{lemma}\label{lemma:L}
Let $a\in\C^*$ and consider the function $L_{a,q}$ defined for
$x\in\C^*\setminus[a;q]$ as above, {\it ie}:
$$
L_{a,q}(x)=\ell_{p}\l(\frac xa\r)=-\frac xa\frac{\theta_p'(-\frac
xa)}{\theta_p(-\frac xa)}\,.
$$
Then, up to an additive constant, $L_{a,q}$ is the only solution
of the $q$-difference equation $y(qx)-y(x)=1$, which is analytic
on $\C^*\setminus[a;q]$ and has only simple poles at $[a;q]$.
\end{lemma}

\begin{proof}
The functional equation for $L_{a,q}$ is obtained deriving the
equation $-\frac xa\theta_p(-p\frac xa)=\theta_p(-\frac xa)$. The
uniqueness (up to a constant) comes from the remark that there are
no non constant elliptic function having only a simple pole in a
fundamental domain.
\end{proof}

\begin{proof}[Proof of Proposition \ref{prop:disc-cont-euler}.]
Let us consider the modular variable change:
$$
q  \longmapsto q^* =e^{4\pi^2/\ln q}\,,\quad x  \longmapsto  x^*
=e^{-2\pi i\frac{\log x}{\ln q}}\,.
$$
In the notation of Lemma \ref{lemma:U} above let
$W(x^*)=U_{\la}(x)$. Then:
$$
W(x^*q^*)-W(x^*)=-2\pi i\frac{q-1}{\ln q},\quad W(x^*e^{-2\pi
i})=W(x^*).
$$
Equivalently, $x^*\mapsto W(x^*)$ represents a uniform solution to
a first order $q^*$-difference equation. By Lemma \ref{lemma:L},
there exists a constant $C\in\C$ such that
$$
W(x^*)=-2\pi i\frac{q-1}{\ln q}L_{(-\la(1-p))^*,q^*}(x^*)+C,\quad
x\in\tilde\Omega_\la\,.
$$
We calculate the constant $C=C_{\la,q}$ setting $x=1-p$ and
$x^*=e^{-2\pi i\frac{\ln(1-p)}{\ln q}}$. Since $e_q(-q/x)$ has a
zero for $x=1-p$, we obtain the exact expression for $C$.
\end{proof}

The main result of this section is:

\begin{thm}\label{thm:moyenne}
For any $x\in\C^*\setminus]-\infty,0[$ we have:
$$
\cE_q(x)=\frac1{\ln q}\int_1^q\cE_q^{[\la]}(x)\frac{d\la}\la.
$$
\end{thm}

The theorem follows from the combination of Corollary
\ref{cor:disc-cont-euler} and the following lemma:

\begin{lemma}
Let $p^*=1/q^*=e^{-\frac{4\pi^2}\ln q}$. For $z$ close enough to
$1$, we have:
$$
\int_1^q\ell_{p^*}\l(\Big(-\frac z\la\Big)^*\r)\frac{d\la}\la=
\int_1^q\ell_{p^*}\l(\Big(-\frac 1\la\Big)^*\r)\frac{d\la}\la\,.
$$
\end{lemma}

\begin{proof}
Let $\mu=-\la^{-1}$. From the identity $x^*=e^{-2\pi i\frac{\log
x}{\ln q}}$ we deduce that:
$$
\frac{d\la}\la=-\frac{d\mu}\mu=-\frac{\ln q}{2\pi
i}\frac{d\mu^*}{\mu^*}\,.
$$
Therefore for $t=\mu^*$ we obtain:
$$
\int_1^q \ell_{p^*}\l(\Big(-\frac z\la\Big)^*\r)\frac{d\la}\la=
\frac{\ln q}{2\pi i}\int_{{\cal
C}_{(-z)^*}}\ell_{p^*}(t)\frac{dt}{t}\,,
$$
where ${\cal C}_{(-z)^*}$ is the positive oriented circle,
centered at 0 and passing through the point $(-z)^*$. Observing
that, for $z$ close enough to $1$, the meromorphic function
$t\mapsto \ell_{p^*}(t)$ has no poles in the annulus between
${\cal C}_{(-z)^*}$ and ${\cal C}_{(-1)^*}$, we conclude applying
Cauchy Theorem.
\end{proof}

\begin{rmk}
In Theorem \ref{thm:moyenne}, we could have replaced the interval
$[1,q]$ with a path of the form $[a,qa]$, for any
$a\in\C\smallsetminus(-\infty,0]$.
\end{rmk}

%%%%%%%%%%%%%%%%%%%%%%%%%%%%%%%%%%%%%%%%%%%%%%%%%%%%%%%%%%%%%%%%%
\subsection{Comparing sums along different spirals}
%%%%%%%%%%%%%%%%%%%%%%%%%%%%%%%%%%%%%%%%%%%%%%%%%%%%%%%%%%%%%%%%%

Notice that if $[\la]=[\mu]$, \ie if $\la$ and $\mu$ are two
complex numbers congruent modulo $q$, the two discrete sums
$\cE_q^{[\la]}$ and $\cE_q^{[\mu]}$ coincide. On the other side,
if $[\la]\not=[\mu]$ these sums are trivially distinct, since the
sets of their poles are distinct. This simple remark underlines a
fundamental difference between the continuous and the discrete
summations. In fact, when we make the direction $d$ vary we are
actually constructing an analytic continuation of $\cE^d_q$ on
$\wtilde\C^\ast$, while when we make $[\la]$ vary, we obtain a
whole family of distinct meromorphic solution of $x^2\dq y+y=x$.
This implies that the ``discrete Stokes phenomenon'' for the
discrete summation has a different nature from the classical
differential Stokes phenomenon. It is described in the following
theorem:

\begin{thm}
For $\lambda,\mu\in\C^*\setminus(1-p)q^{\mathbb Z}$ we have:
$$
\cE_q^{[\la]}-\cE_q^{[\mu]}=\frac{K(\lambda,\mu,x)}{e_p(\frac
qx)}\,,
$$
where:
$$
K(\lambda,\mu,x)=
C\frac{\theta_p(-\frac\lambda\mu)\theta_p(\frac{1-p}x)
\theta_p(\frac{(1-p)\lambda\mu}x)}{\theta_p(\frac{(1-p)\lambda}x)
\theta_p(\frac{(1-p)\mu}x)}\,,
$$
where $C$ is a constant depending only on $q$.
\end{thm}

\begin{proof}
The function $\cE_q^{[\la]}(x)-\cE_q^{[\mu]}(x)$ being solution of
the homogeneous equation $x^2d_qy=-y$, it has the form
$$
\cE_q^{[\la]}-\cE_q^{[\mu]}=\frac{K(\lambda,\mu,x)}{e_p(\frac
qx)}\,,
$$
where $K(\la,\mu,x)$ is $q$-invariant function in each variable
$(x,\la,\mu)$.
\par
We want a more precise description of $K(\la, \mu,x)$. Notice that
$ \cE_q^{[\la]}-\cE_q^{[\mu]}$ has only two spirals of simple
poles: $-(1-p)\la p^{\mathbb Z}$ and $-(1-p)\mu p^{\mathbb Z}$.
Since any $q$-invariant uniform function can be written as a
quotient of Theta functions, we obtain:
$$
K(\la,\mu,x)=\frac{C(\la,\mu)\theta(\frac\alpha
x)\theta(\frac\beta x)}{\theta(\frac{(1-p)\la}
x)\theta(\frac{(1-p)\mu} x)}\,,
$$
where $\alpha\beta=(1-p)^2\lambda\mu$. Moreover the factor
$e_q(-q/x)$ in
$$
K(\la,\mu,x)=e_q(-q/x)(\cE_q^{[\la]}-\cE_q^{[\mu]}),
$$
implies that $K(\la,\mu,x)$ has a spiral of simple zeros at
$(1-p)p^{\mathbb Z}$, which implies that we can chose either $\a$
or $\be$ equal to $-(1-p)$. We conclude that $\{\alpha,
\beta\}=\{-(1-p),-(1-p)\la\mu\}$.
\par
We have to calculate $C(\la,\mu)$. The poles of $K(\la,\mu,x)$
with respect to the variable $\la$ forms two spirals: $-\frac
x{1-p}p^{\mathbb Z}$ and $-p^{\mathbb Z}$, hence:
$$
C(\la,\mu)=\frac{\theta(-\frac\la\mu)}{\theta(\la)}C(\mu)\,.
$$
A similar argument shows that $C(\mu)=C/\theta(p\mu)$.
\end{proof}

\begin{rmk}
One can express the constant $C$ in terms of $q$-series. For
instance, setting $x=1$ and letting $\la\to \mu=1$ in
$K(\la,\mu,x)$, we can express $C$ as a value of a derivative.
\end{rmk}

%%%%%%%%%%%%%%%%%%%%%%%%%%%%%%%%%%%%%%%%%%%%%%%%%%%%%%%%%%%%%%%%%
\subsection{Confluence}\label{subs:confluenceEuler}
%%%%%%%%%%%%%%%%%%%%%%%%%%%%%%%%%%%%%%%%%%%%%%%%%%%%%%%%%%%%%%%%%

\begin{thm}
Let $\cE(x)$ be the sum of the classical Euler series in the
direction ${\mathbb R}^+$. Then $\cE_q(x)\to \cE(x)$ if $q\to 1^+$
for any $x\in\C$ such that $\arg x\in(-\pi,\pi)$ and the
convergence is uniform on the compacts of such a domain.
\end{thm}

\begin{proof}
Notice that for any $t\in{\mathbb R}^+$ we have $e_q(t)\to e^t$
and $e_q(t)\leq e^t$. The dominated convergence theorem applied to
the $q$-Laplace transform in a direction $d\in(-\pi,\pi)$ allows
to conclude. Moreover, the estimate of $e_q(x)$ being uniform with
respect to $d=\arg x$, the uniform convergence on the compacts of
$\{|arg x|<\pi\}$ follows at once.
\end{proof}

\begin{cor}
The same statement holds for $\cE_q^{[\la]}(x)$ when $q\to 1^+$.
\end{cor}

\begin{proof}
The proof results of the combination of Proposition
\ref{prop:disc-cont-euler} relating $\cE_q(x)$ to
$\cE_q^{[\la]}(x)$, Lemma \ref{lemma:log} on the uniform
convergence of the $q$-logarithm, and the theorem above.
\end{proof}

%%%%%%%%%%%%%%%%%%%%%%%%%%%%%%%%%%%%%%%%%%%%%%%%%%%%%%%%%%%%%%%%%%%%%%%%%%%%%%%%%%%%%%%%
%%%%%%%%%%%%%%%%%%%%%%%%%%%%%%%%%%%%%%%%%%%%%%%%%%%%%%%%%%%%%%%%%%%%%%%%%%%%%%%%%%%%%%%%
\section{Generic $q$-Gevrey series}
%%%%%%%%%%%%%%%%%%%%%%%%%%%%%%%%%%%%%%%%%%%%%%%%%%%%%%%%%%%%%%%%%%%%%%%%%%%%%%%%%%%%%%%%
%%%%%%%%%%%%%%%%%%%%%%%%%%%%%%%%%%%%%%%%%%%%%%%%%%%%%%%%%%%%%%%%%%%%%%%%%%%%%%%%%%%%%%%%

We call \emph{generic $q$-Gevrey series} a power series $\hat
f\in\C[[x]]$ satisfying a $q$-difference equation $\Delta \hat
f\in \C\{x\}$ for some analytic $q$-difference operator $\Delta$
of the form:
\begin{equation}\label{equa:pente1a}
\Delta:=a_0(x)(x\sigma_q)^n+a_1(x)(x\sigma_q)^{n-1}+...+a_n(x)\,,
\hbox{~with $a_j\in\C\{x\},\ a_0(0)a_n(0)\not=0$,}
\end{equation}
and $\sgq(f(x))=f(qx)$. This means that the associated Newton
polygon has only one finite slope equal to one (\cf
\cite{ramistoulouse} and \cite{zhangfourier}).
\par
An explicit calculation (\cf also \cite[Lemma 1.1.10]{DV}) shows
that

\begin{lemma}\label{lemma:pente1}
Let $d_q=\frac{\sigma_q-1}{x(q-1)}$ and consider a $q$-difference
operator $\Delta\in\C\{x\}[\sigma_q]$. Then $\Delta$ can be
written as \eqref{equa:pente1a} if and only if it can be rewritten
in the following form:
\begin{equation}\label{equa:pente1b}
\Delta:=b_0(x)(x^2d_q)^n+b_1(x)(x^2d_q)^{n-1}+...+b_n(x), \quad
b_j\in\C\{x\},\ b_0(0)b_n(0)\not=0.
\end{equation}
\end{lemma}

Notice that the $q$-Euler series $\hat \cE(x)$ considered in
previous section is a generic $q$-Gevrey series.

%%%%%%%%%%%%%%%%%%%%%%%%%%%%%%%%%%%%%%%%%%%%%%%%%%%%%%%%%%%%%%%%%%%%%%%%%%%%%%%%%%%%%%%%
\subsection{Two formal $q$-Borel transforms}
\label{subsec:borel}
%%%%%%%%%%%%%%%%%%%%%%%%%%%%%%%%%%%%%%%%%%%%%%%%%%%%%%%%%%%%%%%%%%%%%%%%%%%%%%%%%%%%%%%%

The classical Borel transform associates to each power series
$\sum_{n\ge 0}a_nx^{n+1}$ the \emph{more convergent} (or
\emph{less divergent}) power series $\sum_{n\ge
0}\frac{a_n}{n!}\xi^n$. For the solutions of a $q$-difference
equations, the Gevrey ``scaling factor'' $(n!)^s$ is replaced by
the $q$-Gevrey one: $(q^{sn^2/2})$ (\cf
\cite{Beindex},\cite{ramistoulouse},
\cite{zhangfourier},\cite{zhangthetafunction}). Indeed, in the
literature there are (at least) two $q$-analogs of the factorial
$n!$, namely $[n]_q^!$ and $q^{n(n-1)/2}$. The reason for this
dichotomy is the identity
$$
[n]_q^! =\frac{(q;q)_n}{(1-q)^n}
=\frac{(p;p)_n}{(1-p)^n}q^{n(n-1)/2} =[n]_p^!q^{n(n-1)/2}\,,
$$
which implies that
\begin{equation}\label{equa:equiv}
[n]_q^!\sim\frac{q^{n(n-1)/2}}{(1-p)^n}\,, \hbox{ when
}n\to+\infty\,.
\end{equation}

Let us consider the following two formal Borel transforms,
associated to those $q$-factorials:
$$
{\cal B}_q\ :\quad x\C\l[[x]\r]\to \C\l[[\xi]\r],\ \sum_{n\ge
0}a_n x^{n+1}\longmapsto \sum_{n\ge 0}\frac{a_n}{[n]_q^!}\xi^n\,
$$
and
$$
B_q\ :\quad x\C\l[[x]\r]\to \C\l[[\xi]\r],\ \sum_{n\ge 0}a_n
x^{n+1}\longmapsto \sum_{n\ge 0}\frac{a_n}{q^{n(n-1)/2}}\xi^n.
$$
Following J.P. Ramis \cite{ramistoulouse} we set:

\begin{defn}
An entire function $\phi$ is said to have \emph{$q$-exponential
growth of order 1 at $\infty$} if there exist two constants
$K>0$ and $\mu>0$ such that
$$
\vert \phi(x)\vert<K\vert x\vert ^{\mu}e_q(\vert x\vert), \quad
\vert x\vert\to\infty.
$$
\end{defn}

Remark that the function $e_q(\vert x\vert)$ can be replaced by
$e^{\ln^2\vert x\vert/(2\ln q)}$.

\begin{lemma}[{\cite[Prop. 2.1]{ramistoulouse}}]
\label{lemma:croissanceqexp} Let $\E_q$ be the set of all the
entire functions having a $q$-exponential growth of order 1 at
$\infty$, and let $\C\{x\}$ be the set of all power series having
a positive convergence radius. Then $\E_q={\cal
B}_q(x\C\{x\})=B_q(x\C\{x\})$.
\end{lemma}

The following function space $\H_q$ has been introduced in
\cite{zhanggroningen} and \cite{RZ}; see also \cite{zhangfourier}.

\begin{defn} %(\cite{zhanggroningen}, \cite{RZ})
For any $\lambda\in\C^*$, let $[\lambda;q]=\lambda q^{\mathbb Z}$.
\begin{enumerate}

\item
A germ of  function $\phi$ analytic at $0$ is said to belong to
$\H^{[\lambda;q]}$ if there exist a domain $\Omega\subset\C$ and a
real number $r>0$ such that:

\begin{itemize}
\item $\cup_{m\ge 0}\{x\in\C: \vert x-\lambda q^m\vert<rq^m\}\subset\Omega$;
\item $\phi$ can be continued to be an analytic function on
$\Omega$ with a $q$-exponential growth of order 1 at infinity.
\end{itemize}

\item
A germ of  function $\phi$ analytic at $0$ is said to belong to
$\H_q$ if there exist a finite set $\Lambda\subset\C^*$ such that
$\phi\in\H^{[\lambda;q]}$ for any
$\lambda\in\C^*\setminus\Lambda$.
\end{enumerate}

\end{defn}

\begin{prop}
Let $\hat f\in x\C[[x]]$, $\lambda\in\C^*$ and let $\psi=\cB_q\hat
f$, $\phi=B_q\hat f$. Then $\psi\in\H^{[\lambda;q]}$ if and only
if $\phi\in\H^{[(1-p)\lambda;q]}$.
\par
Moreover, the map $ \sum_{n\ge 0}a_nx^n\mapsto\sum_{n\ge
0}a_n[n]_p^!x^n $ induces an automorphism of the vector space
$\H_q$.
\end{prop}

\begin{proof}
Use the integral representation for the corresponding Hadamard
product.
\end{proof}

\begin{defn}
Let $\lambda\in\C^*$ and let $\hat f\in x\C[[x]]$.
\begin{enumerate}

\item
We set $\C\{x\}_q^{[\lambda;q]}=\cB_q^{-1}\H^{[\lambda/(1-p);q]}=
B_q^{-1}\H^{[\lambda;q]}$. We say that $\hat
f\in\C\{x\}_q^{[\lambda;q]}$ is \emph{$q$-Borel summable along
$[\lambda;q]$}.

\item
Each element of $\DS\hat f:=\{[\lambda;q]:\hat
f\notin\C\{x\}_q^{[\lambda;q]}\}$ will be called \emph{singular
direction of $\hat f$}.

\item
If $\DS\hat f $ is finite, $\hat f$ is called \emph{$q$-Borel
summable} and we denote $\C\{x\}_q$ the set of $q$-Borel summable
series.
\end{enumerate}
\end{defn}

\begin{thm} [{\cite[Thm. 1.2.1]{zhanggroningen}}]
Every generic $q$-Gevrey series is $q$-Borel summable.
\end{thm}

\subsection{Different kinds of $q$-exponential summation}
\label{subsec:summations}

The classical Borel-Laplace exponential summation is based on the
Euler's integral representation of the Gamma function, namely
\begin{equation}\label{equa:EulerGamma}
n!=\int_0^\infty e^{-t}t^{n+1}\frac{dt}t\,.
\end{equation}
In the definition of a $q$-summation procedure one must be guided
by the $q$-analogs of this last integral, question investigated
since Jackson, Wigeret, Watson, etc...

We recall the following $q$-analogs of the integral representation
of the Euler Gamma function.

\begin{prop}\label{prop:qfactorial}
For $d\in(-\pi,\pi)$ and $\lambda\notin(-q^{\mathbb Z})$ we have:
\begin{equation}\label{equa:qfactorial}
    [n]_q^!=\frac{q-1}{\ln q}\, \int_0^{\infty e^{id}}\frac{t^n}{e_q(qt)}\,dt
    =\frac{q-1}{\ln q}\int_0^{\infty e^{id}}{t^n}{e_p(-qt)}\,dt\,,
\end{equation}

\begin{equation}\label{equa:qfactorialdq}
    [n]_q^!=q\int_{\la p^{\mathbb Z}}\frac{t^n}{e_q(qt)}\,d_pt
    =q\int_{\la p^{\mathbb Z}}t^n e_p(-qt)\,d_pt\,,
\end{equation}

\begin{equation}\label{equa:qfactorialq}
q^{n(n-1)/2}=\frac{q}{\ln q}\int_0^{\infty
e^{id}}\frac{t^n}{\theta_p(qt)}\,dt\,,
\end{equation}

\begin{equation}\label{equa:qfactorialqq}
q^{n(n-1)/2}=\frac{q}{1-p}\int_{\la p^{\mathbb
Z}}\frac{t^n}{\theta_p(qt)}\,d_pt\,.
\end{equation}

\end{prop}

\begin{proof}
For the proof of the identities above \cf \cite[pages
549-550]{AA}. More precisely, letting $c\to n+1$, $b\to 1$ and
$a\to 0$ (resp. letting $c\to n+1$, $b,a\to 0$) in
$$
\int_0^\infty
x^{c-1}\frac{(-ax;p)_\infty(-pb/x;p)_\infty}{(-x;p)_\infty(-p/x;p)_\infty}dx=
\frac{(ab;p)_\infty(p^c;p)_\infty(p^{1-c};p)_\infty\pi}
{(bp^c;p)_\infty(ap^{-c};p)_\infty(p;p)_\infty\sin(\pi c)},
$$
one gets the formula
$$
\int_0^\infty \frac{x^{n}}{(-x;p)_\infty}dx=(\ln q)\ (p;p)_n
p^{-n(n+1)/2}$$
$$
\l(\hbox{resp.~} \int_0^\infty
\frac{x^{n}}{(-x;p)_\infty(-p/x;p)_\infty}dx=(\ln q)\
p^{-n(n+1)/2}\r),
$$
which yields \eqref{equa:qfactorial} (resp.
\eqref{equa:qfactorialq}). Similarly, the formulae
\eqref{equa:qfactorialdq}  and  \eqref{equa:qfactorialqq} can be
viewed as special cases of
$$
\int_0^\infty x^{c-1}
\frac{(-ax;p)_\infty(-bp/x;p)_\infty}{(-x;p)_\infty(-p/x;p)_\infty}d_px=
\frac{(1-p)(p;p)_\infty(-p^c;p)_\infty(-p^{1-c};p)_\infty(ab;p)_\infty}
        {(-1;p)_\infty(-p;p)_\infty(ap^{-c};p)_\infty(bp^{c};p)_\infty}\,.
$$
\end{proof}

\begin{rmk}
In particular, \eqref{equa:qfactorialqq} and
\eqref{equa:qfactorialq} have been studied in \cite{RZ} and
\cite{zhanggroningen} as starting points for the corresponding
summation procedures. Other kinds of $q$-summation are considered
in \cite{zhangfourier} and \cite{MarotteZhang}.
\end{rmk}

Let $d\in[-\pi,\pi)$. We will identify $d$ to the half line
$[0,\infty e^{id}):={\mathbb R}^+e^{id}$.

\begin{defn}
We set $\H_q^d=\cap_{\lambda\in (0,\infty
e^{id})}\H^{[\lambda;q]}$.
\end{defn}

\begin{rmk}
The functional space $\H_q^d$ is exactly the space $\H_{q;1}^d$
introduced in \cite{zhangfourier}.
\end{rmk}

Let $\lambda\in\C^*$ and $d\in[-\pi,\pi)$. According to
Proposition \ref{prop:qfactorial}, the following four
\emph{$q$-Laplace transforms} are well defined:
$$
\forall \phi\in\H^{[\lambda;q]},\quad L_q^{[\lambda]}\phi=\frac
q{1-p}\int_{\la p^{\mathbb Z}}\frac{\phi(\xi)} {\theta_p(q\frac
\xi x)}\,{d_p\xi}\,;
$$
$$
\forall \phi\in\H^{[\lambda/(1-p);q]},\quad
 \cL_q^{[\lambda]}\phi=\frac q{1-p}
    \int_{\la p^{\mathbb Z}}\frac{\phi(\xi/(1-p))}{e_q(q\frac \xi {(1-p)x})}\,{d_p\xi};
$$
 $$
\forall \phi\in\H^d_q,\quad L_q^d\phi=\frac q{\ln
q}\int_0^{e^{id}\infty}\frac{\phi(\xi)}
    {\theta_p(q\frac \xi x)}d\xi\,,\quad \cL_q^{d}\phi=\frac{q-1}{\ln q}
    \int_0^{e^{id}\infty}\frac{\phi(\xi)}{e_q(q\frac \xi x)}d\xi\,.
$$

\begin{defn}\label{defn:summations}
\begin{enumerate}

\item
If $\hat f\in\C\{x\}_q^{[\lambda;q]}$, we define its sums in the
direction $[\lambda;q]$ as follows:
$$
{\cal S}_q^{[\lambda]}\hat f=\cL_q^{[\lambda]}(\cB_q\hat f),\quad
S_q^{[\lambda]}\hat f=L_q^{[\lambda]}(B_q\hat f).
$$

\item
If $\hat f\in\C\{x\}_q^{d}$, define its sums in the direction $d$
as follows:
$$
{\cal S}_q^{d}\hat f=\cL_q^{d}(\cB_q\hat f),\quad S_q^{d}\hat
f=L_q^d(B_q\hat f).
$$
\end{enumerate}

\end{defn}

\begin{rmk}~\label{rmk:moyenne}
\begin{trivlist}
\item $\bullet$
The summation procedures $\hat f\to S^{[\lambda;q]}\hat f$ and
$\hat f\to S^{d}\hat f$ are introduced in \cite{zhanggroningen}
and \cite{zhangthetafunction}: they have many good asymptotic
properties.

\item $\bullet$
Suppose that $\hat f$ is $q$-summable and that  $d$ is not a
singular direction. Then we have the following \emph{formal}
equality (meaning that we exchange carelessly the infinite sum and
the integral):
\begin{equation}\label{eq:moyenneq1}
\cS_q^d\hat f=\frac 1{\ln q}
\int_{e^{id}}^{qe^{id}}\cS_q^{[\la]}\hat f\frac{d\la}\la
\end{equation}
and
\begin{equation}\label{eq:moyenneq2}
S_q^d\hat f=\frac 1{\ln q}
\int_{e^{id}}^{qe^{id}}S_q^{[\la]}\hat f\frac{d\la}\la\,.
\end{equation}
To prove that this identity is not only formal, but analytic, one
would like to apply the dominated convergence theorem:
unfortunately the dominated convergence is a little bit delicate
for a general $\hat f$, since we don't really control the spirals
of poles of the discrete $q$-Borel sums. Anyway, we will prove
\eqref{eq:moyenneq1} and \eqref{eq:moyenneq2} for a generic
$q$-Gevrey series (\cf Theorem \ref{thm:sums}).
\end{trivlist}
\end{rmk}

At this stage a natural question arises:
$$
\hbox{Do we have $S^{[\lambda;q]}\hat f ={\cal
S}^{[\lambda;q]}\hat f$ and $S^d\hat f={\cal S}^d\hat f$?}
$$
The answer is clear, and trivially positive, if $\hat f(x)$ is a
germ of analytic function at zero: in this case all the sums of
$\hat f(x)$ coincide with $f$.

The rest of the paper is devoted to the proof of the following
theorem:

\begin{thm}\label{thm:sums}
Let $\hat f$ be a generic $q$-Gevrey series and let
$\lambda\in\C^*$, $d\in[-\pi,\pi)$. Assume that
$\lambda\notin\DS(\hat f)$ and $(0,e^{id}\infty)\cap\DS(\hat
f)=\emptyset$. Then
\begin{enumerate}
\item
$S_q^{[\la]}\hat f=\cS_q^{[\la]}\hat f$, on a convenient domain
$\Omega$.

\item
$S_q^d\hat f=\cS_q^d\hat f$ on a convenient sector containing the
direction $d$.
\end{enumerate}
Moreover we have:
$$
\cS^d_q\hat f=\frac1{\ln q}\int_{e^{id}}^{qe^{id}}
\cS^{[\la]}_q\hat f\frac{d\la}\la\,.
$$
\end{thm}

\begin{rmk}
Theorems \ref{thm:identity} and \ref{thm:moyenne} are a special
case of Theorem \ref{thm:sums}.
\end{rmk}

Before giving a proof of Theorem \ref{thm:sums} in
\S\ref{subs:proof}, we make a digression about two essential
ingredients of the proof: first we prove the theorem in the
special case of the Tschakaloff series; then we introduce a
functional space that allows to read, in certain sense, any
$q$-Gevrey series as a finite linear combination of some modified
Tschakaloff series.

%%%%%%%%%%%%%%%%%%%%%%%%%%%%%%%%%%%%%%%%%%%%%%%%%%%%%%%%%%%%%%%%%%%%%%%%%%%%%%%%%%%%%%%%
\subsection{The Tschakaloff series}
%%%%%%%%%%%%%%%%%%%%%%%%%%%%%%%%%%%%%%%%%%%%%%%%%%%%%%%%%%%%%%%%%%%%%%%%%%%%%%%%%%%%%%%%

Let us consider another $q$-analogue of the Euler series:
$$
T_q(x)=\sum_{n\geq 0}q^{n(n-1)/2}x^{n+1}\,,
$$
called the Tschakaloff series or the partial Theta function. It
satisfies the $q$-difference equation
$$
xT_q(qx)-qT_q(x)=-qx\,,
$$
that can also be rewritten in the form $x^2(q-1)\dq y+(x-q)y=-qx$.
\par
The Borel transforms of $T_q$ are:
$$
\psi(\xi)=\cB_q(T_q)=\sum_{n\geq
0}\frac{\xi^n}{[n]_p^!}=e_p(\xi)=\l((1-p)\xi;p\r)^{-1}
$$
and
$$
\phi(\xi)=B_q(T_q)=\frac 1{1-\xi}\,.
$$

\begin{prop}\label{prop:Tschakaloff}
Let us fix $\la\notin[-1;q]$. Then
$\cS_q^{[\la]}T_q=S_q^{[\la]}T_q$.
\end{prop}

\begin{proof}
The definition of the Jackson integral (\cf
\S\ref{sec:jacobsonintegral}) and the Jacobi triple product
formula (\cf \eqref{eq:Jacobi}), plus the development of the
$q$-exponential $e_q(x)$ as an infinite product (\cf
\eqref{eq:q-exp}), imply that:
$$
\begin{array}{rcl}
\cL_q^{[\la]}\psi &=&\ds\frac q{1-p}\int_{\la p^{\mathbb
Z}}\frac{e_p(\xi)}{e_q\l(\frac{q\xi}{(1-p)x}\r)}d_p\xi
%\\&=&
=\ds\frac q{1-p}\int_{\la p^{\mathbb Z}}\l(\xi,-\frac{q\xi}{x};p\r)_\infty^{-1}d_p\xi\\ \\
%&=&\ds\frac q{1-p}(1-p)\la\sum_{n\in\Z}p^n\l(p^n\la,-\frac{p^{n-1}\la}{x};p\r)_\infty^{-1}\\
%&=&\ds \la\sum_{n\in\Z}p^{n-1}\l(p^n\la,-\frac{p^{n-1}\la}{x};p\r)_\infty^{-1}\\
&=&\ds \la(p;p)_\infty\sum_{n\in{\mathbb Z}}
    \frac{p^n(-p^{1-n}x /\la;p)_\infty}{(p^{n+1}\la;p)_\infty\theta_p(p^n\la/x )}\,.
\end{array}
$$
Since $\theta_p(x)=p^{n(n-1)/2}x^n\theta_p(p^nx)$ for any
$n\in{\mathbb Z}$, we obtain:
$$
\cS_q^{[\la]}T_q=\cL_q^{[\la]}\psi
=\frac{\la(p;p)_\infty}{\theta_p(\la/x )}\sum_{n\in{\mathbb Z}}
\frac{\l(-p^{1-n}\frac x \la;p\r)_\infty}{(p^{n+1}\la;p)_\infty}
p^{n(n+1)/2}\l(\frac \la x\r)^n\,.
$$
On the other hand we have:
$$
\begin{array}{rcl}
S_q^{[\la]}T_q=L_q^{[\la]}\phi
&=&\ds \frac q{1-p}\int_{p^{\mathbb Z}\la}\frac{d_p\xi}{(1-\xi)\theta_p(q\xi/x)}\\ \\
&=&\ds \la\sum_{n\in{\mathbb Z}}\frac{p^n}{(1-p^{n+1}\la)\theta_p(p^n\la/x)}\\ \\
&=&\ds \frac\la{\theta_p(\la/x)}\sum_{n\in{\mathbb
Z}}\frac{p^{n(n+1)/2}}{1-p^{n+1}\la}\l(\frac \la x\r)^n
\end{array}
$$
A straightforward calculation of the residues of
$(x,p)_\infty^{-1}$ at $x=p^{-k}$, for $k\geq 0$, gives the
formula:
$$
\frac1{(x;p)_\infty}=\frac1{(p;p)_\infty} \sum_{k\ge
0}\frac{a_k}{1-xp^k}\,, \hbox{ with } a_k=\frac1{(p^{-k};p)_k}=
\frac{(-1)^kp^{k(k+1)/2}}{(p;p)_k}\,.
$$
Therefore we obtain:
$$
\begin{array}{rcl}
\cS_q^{[\la]}T_q &=&\ds \frac\la{\theta(\frac\la
x)}\sum_{\ell\in{\mathbb Z}}\frac{p^{\ell(\ell+1)/2}}{1-p^\ell
\la} \l(\frac\la x\r)^\ell (a_\ell;p)_\infty\
{}_0\phi_1(-;a_\ell;p,a_\ell),\quad a_\ell=-p^{1-\ell}\frac
x\la\,.
\end{array}
$$
The Ramanujan formula (\cf \cite[Thm. 4.4]{zh2}):
$$
(x;p)_\infty\ {}_0\phi_1(-;x;p,x)=1\,, \hbox{ for any $x\notin
q^{-{\mathbb N}}$,}
$$
implies that $\cL_q^{[\la]}\psi=L_q^{[\la]}\phi$.
\end{proof}

%%%%%%%%%%%%%%%%%%%%%%%%%%%%%%%%%%%%%%%%%%%%%%%%%%%%%%%%%%%%%%%%%%%%%%%%%%%%%%%%%%%%%%%%
\subsection{The functional space $H$ in the Borel plane}
\label{subs:spaceH}
%%%%%%%%%%%%%%%%%%%%%%%%%%%%%%%%%%%%%%%%%%%%%%%%%%%%%%%%%%%%%%%%%%%%%%%%%%%%%%%%%%%%%%%%

We recall that the $q$-Borel transform $B_q$ associates to a
power series $\hat f=\sum_{n\ge 0}a_nx^{n+1}\in\C[[x]]$ the
power series $\phi=\sum_{n\ge
0}{a_n}q^{-n(n-1)/2}\xi^n\in\C[[\xi]]$. As we have already pointed
out, the $q$-Borel transform of a generic $q$-Gevrey series admits
a positif radius of convergence and can be continued to an
analytic function in the whole complex plane minus a finite number
of sets of the form $\la q^{\mathbb N}$ (\cf \cite{zhangfourier}
and \cite{zhanggroningen}).
\par
In this section we want to prove that every generic $q$-Gevrey
series can be expressed by means of ``modified Tschakaloff
series''. Our strategy consists in proving that the $q$-Borel
transform of any generic $q$-Gevrey series admits an elementary
decomposition, by studying the $q$-convolution product of suitable
entire functions by a rational functions. This leads to the
construction of a functional space which is somehow spanned by the
$q$-Borel transforms of the modified Tschakaloff series.

\begin{defn}
We call $q$-convolution product the following bilinear operator:
$$
\begin{array}{rccc}
*_q:& \C\{\xi\}\times\C\{\xi\} & \longrightarrow & \xi\C\{\xi\}\\ \\
&\xi^n*_q\xi^m & \longmapsto & q^{-(nm+n+m+1)}\xi^{n+m+1}
\end{array}\,.
$$
\end{defn}

A direct calculation shows that
\begin{trivlist}

\item 1.
If $\phi=\sum_{n\ge 0}\phi_n\xi^n\in\C\{\xi\}$ and
$\psi\in\C\{\xi\}$, then (\cite[1.4.3, where
$s=1$)]{MarotteZhang})
$$
\phi*_q\psi(\xi)=\sum_{\ge
0}\phi_0q^{-n-1}\xi^{n+1}\psi(q^{-n-1}\xi).
$$

\item 2.
$B_q(\hat f\hat g)=B_q(\hat f)*_qB_q(\hat g)$.

\end{trivlist}

Let $K$ be the set of rational functions bounded at zero and let
${\mathbb E}_q $ be the set of all entire functions admitting at
most a $q$-exponential growth of order 1 at the infinity. We know
that $K\cap {\mathbb E}_q =\C[\xi]$ and ${\mathbb E}_q
=B_q(x\C\{x\})$ (\cf [Ram92]). Notice that the formula $B_q(\hat
f\hat g)=B_q(\hat f)*_qB_q(\hat g)$ identifies $({\mathbb E}_q
,*_q)$ to a commutatif sub-ring of $(\C\{\xi\},*_q)$.

\begin{defn}
We define the functional space $H:=\cup_{n\ge 0}H_n$ in the
following way:
$$
H_{-1}=\{1\},\quad H_0=K,\quad H_1={\mathbb E}_q *_qK:=\{\phi*_qr:
\phi\in {\mathbb E}_q ,r\in K\}
$$
and, for any integer $n\ge 1$,
$$
H_{2n}=KH_{2n-1}:=\{ru:r\in K,u\in H_{2n-1}\},\quad
H_{2n+1}={\mathbb E}_q *_qH_{2n}:=\{\phi*_qu:\phi\in {\mathbb E}_q
,u\in H_{2n}\}.
$$
\end{defn}

\begin{prop}
For any $(r,\phi,u)\in K\times {\mathbb E}_q \times H$, we have
$(ru,\phi*_qu)\in H\times H$. In other words, the functional space
$H$ is a $(K,{\mathbb E}_q )$-bimodule.
\end{prop}

\begin{proof}
It follows immediately from the definition of $H$. Indeed,  if
$n\le m$, then $H_n\subset H_m$. So, we can suppose that
$(r,\phi,u)\in K\times {\mathbb E}_q \times H_n$ and hence,
$(ru,\phi*_qu)\in H_{n+2}\times H_{n+2}\subset H\times H$.
\end{proof}

\begin{thm}\label{thm:decomposition}
For any $u\in H$, there exist $\phi_0$, $\phi_1,\dots,\phi_n\in
{\mathbb E}_q $ and $r_0$, $r_1,\dots,r_n\in K$ such that
$$
u=\phi_0+r_0+\phi_1*_qr_1+...+\phi_n*_qr_n.
$$
Moreover, we can suppose that $r_1$, ..., $r_n$ are rational
functions of the form $\frac1{(\xi-\la_i)^{\nu_i}}$, where $c_i$,
$\la_i\in\C^*$ and $\nu_i\in{\mathbb N}$.
\end{thm}

\begin{proof}
Since $H_m\subset H_{m+1}$, for any $u\in H$ there exists
$m\in{\mathbb N}$ such that $u\in H_m$. So we can prove the
theorem by induction on $m$. The cases $m=0$ and $m=1$ are
trivial. Suppose that $u\in H_{m+1}$. Then there exists
$(r,v,\phi)\in K\times H_m\times {\mathbb E}_q $ such that on of
the following two cases occurs:
\begin{trivlist}
\item (1) $u=rv$,
\item (2) $u=\phi*v$,
\end{trivlist}
and, by inductional hypothesis,
$v=\phi_0+r_0+\sum_{j=1}^m\phi_j*_qr_j$. The proof in the case (2)
is straightforward, since
$\phi*_q(\phi_j*_qr_j)=(\phi*_q\phi_j)*_qr_j$. In the case (1), we
need the following elementary lemma.

\begin{lemma}
For any $(a,b,\ell,n)\in\C^*\times\C^*\times{\mathbb
N}\times{\mathbb N}$ such that $a\not=b$, the following
decomposition holds:
$$
\frac1{(x-a)^\ell}\frac1{(x-b)^n}=\sum_{k=0}^{\ell-1}\frac{(n+k)!}
{(a-b)^{n+k}k!}\frac{(-1)^k}{(x-a)^{\ell-k}}+\sum_{k=0}^{n-1}
\frac{(\ell+k)!}{(b-a)^{\ell+k}k!}\frac{(-1)^k}{(x-b)^{n-k}}.
$$
\end{lemma}

\begin{proof}
It enough to take the $(\ell-1)$-th derivative with respect to $a$
and the $n-1$-th derivative with respect to $b$ in the formula
$$
\frac1{x-a}\frac1{x-b}=\frac1{a-b}\l(\frac1{x-a}-\frac1{x-b}\r).
$$.
\end{proof}

Let us go back to the proof of the Theorem \ref{thm:decomposition}.
By linearity, it enough to consider a product of the form
$r(\phi*_qr')$, with $r=\frac 1{(\xi-\la)^\ell}$, $r'=\frac
1{(\xi-\mu)^n)}$ and $\phi=\sum_{k\ge 0}\phi_k\xi^k\in {\mathbb
E}_q $. Since
$$
r'*_q\phi(\xi)=\sum_{k\ge
0}\phi_kq^{-k-1}\frac{\xi^{k+1}}{(q^{-k-1}\xi-\mu)^n}\,,
$$
the decomposition follows from the lemma above.
\end{proof}

\begin{cor}
For any $u\in H$, there exist $\la_1,\dots,\la_n\in\C^*$ such that
$u$ is analytic on the domain
$\C\setminus(\cup_{i=1}^m\la_iq^{\mathbb N})$ and the function $U$
defined by
$$
U(\xi)=u(\xi)\prod_{i=1}^n\prod_{m\ge 0}\l(1-\frac\xi{\la_iq^m}\r)
$$
can be continued to an entire function that has at most a
$q$-exponential growth of order $n+1$ at the infinity.
\end{cor}

The corollary results from the combination of the theorem above
and the following lemma:

\begin{lemma}
Let $\phi\in {\mathbb E}_q $, $r=\frac1{(\xi-\la)^n}$, $n\ge 1$
and $\la\in\C^*$. Then $\phi*_qr$ admits  $\la q^{\mathbb N}$ as
set of poles and there exist $C>0$, $m>0$ such that, for any
$\epsilon>0$,
$$
\vert \xi q^{-n}-\la\vert>\epsilon\quad \Longrightarrow\quad \vert
\phi*_qr(\xi)\vert <\frac C{\epsilon^n}\vert\xi\vert^m\vert
e^{\frac{(\log x)^2}{2\ln q}}\vert.
$$
\end{lemma}

\begin{proof}
Let $\phi=\sum_{k\ge 0}\phi_n\xi^k$. Since $\phi\in {\mathbb E}_q
$, there exist $A$, $B>0$ such that
$$
\forall k\in{\mathbb N},\quad \vert\phi_k\vert<AB^nq^{-n(n+1)/2}.
$$
On the other hand,
$$\phi*_qr=\sum_{k\ge 0}\phi_kq^{-k-1}\xi^{k+1}r(q^{-k-1}\xi),
$$
which implies directly the lemma.
\end{proof}

%%%%%%%%%%%%%%%%%%%%%%%%%%%%%%%%%%%%%%%%%%%%%%%%%%%%%%%%%%%%%%%%%%%%%%%%
\subsection{Proof of Theorem \ref{thm:sums}}\label{subs:proof}
%%%%%%%%%%%%%%%%%%%%%%%%%%%%%%%%%%%%%%%%%%%%%%%%%%%%%%%%%%%%%%%%%%%%%%%%

We start by proving the following preparatory result:

\begin{prop}
If $\hat f\in x\C[[x]]$ is a generic $q$-Gevrey series, then its
$q$-Borel transform belongs to $H$.
\end{prop}

\begin{proof}
Let $\Delta$ be a linear analytic $q$-difference operator such
that $\Delta\hat f=g\in x\C\{x\}$. We know that $\Delta$ admits an
analytic factorization (\cf \cite[Prop. 5.1.4]{zhangfourier},
\cite[Thm. 1.2.1]{sauloyfiltration}):
\begin{equation}\label{equa:factor}
\Delta=(x\sigma_q-\la_1)h_1(x\sigma_q-\la_2)h_2...(x\sigma_q-\la_n)h_n,\quad
\la_j\in\C,\ h_j\in\C\{x\},\ h_j(0)=1\,.
\end{equation}
We suppose that we have chosen $n$ minimal and let us prove the
statement by induction on $n$. We consider first of all the case
$n=1$: we suppose that $(x\sgq -\la_1)h_1\hat f=g$, with
$B_q(g)\in H$\footnote{Notice that we are assuming that $B_q(g)\in
H$ and not $B_q(g)\in {\mathbb E}_q $.}. This implies that
$B_q(h_1\hat f)\in H$, since $B_q((x\sigma_q-\la_1)h_1\hat f)
=(q\xi-\la_1)B_q(h_1\hat f)$. Therefore
$$
B_q(xh_1(qx)\hat f(qx))=B_q(g)-B_q(\la_1h_1\hat f)\in H\,,
$$
with $B_q(xh_1(qx))\in {\mathbb E}_q $ and $xh_1(qx)\hat f(qx)\in
x^2\C[[x]]$. So $x\hat f(qx)=\wtilde g h_1(qx)^{-1}$ and
$B_q(x\hat f(qx))=B_q(\wtilde g/x)\ast_q B_q(xh_1(qx)^{-1})\in H$.
Finally $B_q(\hat f)\in H$.
\par
For $n>1$, the inductive hypothesis implies that
$B_q((x\sigma_q-\la_n)h_n\hat f)\in H$, and hence that $B_q(\hat
f)\in H$.
\end{proof}

\begin{proof}[Proof of Theorem \ref{thm:sums}]
Applying Theorem \ref{thm:decomposition} to $B_q(\hat f)$, we can
write $\hat f$ as follows:
$$
\hat f=f_0+\hat e_0+f_1\hat e_1+...+f_n\hat e_n,
$$
where $f_0,\dots,f_n\in x\C\{x\}$, $B_q(\hat e_0)\in K$ and, for
$i=1$,..., $n$, $B_q(\hat e_i)=\frac1{(\xi-\la_i)^{\nu_i}}$. So it
is enough to prove Theorem \ref{thm:sums} for a modified
Tschakaloff series (\cf \cite[Prop. 1.4.2]{MarotteZhang}), \ie
under the assumption $B_q(\hat f)=\frac1{(\xi-\la)^n}$. By
replacing $\xi$ by $-\la\xi$, we can suppose that $\la=-1$. The
case $n=1$ corresponds  exactly to the Tschakaloff divergent
series $\hat \cE$, and the result is stated  in Proposition
\ref{prop:Tschakaloff}. If $n>1$, by considering $\hat
f(x)=\frac{(-1)^{n-1}}{(n-1)!x^{n-1}}{\partial^{n-1}\ }{\partial
a^{n-1}}\hat \cE(ax) \Vert_{a=1}$, we can easily deduce the wanted
result by the help of the dominated convergence theorem.
\par
Concerning the second part of the statement of Theorem
\ref{thm:sums}, the decomposition above allows once again to
reduce to the case of the Tschakaloff series. The dominated
convergence theorem applies with no difficulties to this explicit
case (\cf Remark \ref{rmk:moyenne}).
\end{proof}

%\makeatletter
%\renewcommand{\thesection}{\@Alph\c@section}
%\makeatother
%\setcounter{section}{2}
%\section{}

%%%%%%%%%%%%%%%%%%%%%%%%%%%%%%%%%%%%%%%%%%%%%%%%%%%%%%%%%%%%%%%%%%%%%%%%%%%%%%%%%%%%%%%%
%\nocite{gazette,RamisDivergentes,raminet,lodayexpositiones,
%  lodaygazette,malgrangeexpositiones,Beindex}
%\bibliography{qdiff}

\begin{thebibliography}{DVRSZ03}

\bibitem[AAR99]{AA}
G.~E.~Andrew, R.~Askey, and R.~Roy.
\newblock {\em Special functions}, volume~71 of {\em Encyclopedia of
  Mathematics and Its Applications}.
\newblock Cambridge University Press, Cambridge, 1999.

\bibitem[And00]{andreannalsII}
Y.~Andr\'{e}.
\newblock S\'eries {G}evrey de type arithm\'etique. {I}. {T}h\'eor\`emes de
  puret\'e et de dualit\'e.
\newblock {\em Ann. of Math. (2)}, 151(2):705--740, 2000.

\bibitem[Bar08]{Barnes}
E.W.~Barnes.
\newblock A new development of the theory of the hypergeometric functions.
\newblock {\em Proceedings of the London Mathematical Society}, 6(2):141--177,
  1908.

\bibitem[B\'{e}z92]{Beindex}
J.-P.~B\'{e}zivin.
\newblock Sur les \'equations fonctionnelles aux {$q$}-diff\'erences.
\newblock {\em Aequationes Mathematicae}, 43(2-3):159--176, 1992.

\bibitem[DSK05]{DeSoleKac}
A.~De~Sole and V.~G.~Kac.
\newblock On integral representations of {$q$}-gamma and {$q$}-beta functions.
\newblock {\em Atti della Accademia Nazionale dei Lincei. Classe di Scienze
  Fisiche, Matematiche e Naturali. Rendiconti Lincei. Serie IX. Matematica e
  Applicazioni}, 16(1):11--29, 2005.

\bibitem[DV02]{DV}
L.~Di~Vizio.
\newblock Arithmetic theory of {$q$}-difference equations. {T}he {$q$}-analogue
  of {G}rothendieck-{K}atz's conjecture on {$p$}-curvatures.
\newblock {\em Inventiones Mathematicae}, 150(3):517--578, 2002.

\bibitem[DVRSZ03]{gazette}
L.~Di~Vizio, J.-P.~Ramis, J.~Sauloy, and C.~Zhang.
\newblock \'{E}quations aux {$q$}-diff\'erences.
\newblock {\em Gaz. Math.}, (96):20--49, 2003.

\bibitem[EMOT81]{erd}
A.~Erd\'{e}lyi, W.~Magnus, F.~Oberhettinger, and F.~G.~Tricomi, editors.
\newblock {\em Higher transcendental functions. {V}ol. {I}}.
\newblock Robert E. Krieger Publishing Co. Inc., Melbourne, Fla., 1981.
\newblock Based on notes left by Harry Bateman. Reprint of the 1953 original.

\bibitem[GR90]{GR}
G.~Gasper and M.~Rahman.
\newblock {\em Basic hypergeometric series}, volume~35 of {\em Encyclopedia of
  Mathematics and its Applications}.
\newblock Cambridge University Press, Cambridge, 1990.
\newblock With a foreword by Richard Askey.

\bibitem[LR90]{lodaygazette}
M.~Loday-Richaud.
\newblock Introduction \`{a}\ la multisommabilit\'{e}.
\newblock {\em Gazette des Math\'ematiciens}, (44):41--63, 1990.

\bibitem[LR95]{lodayexpositiones}
M.~Loday-Richaud.
\newblock Solutions formelles des syst\`emes diff\'erentiels lin\'eaires
  m\'eromorphes et sommation.
\newblock {\em Expositiones Mathematicae. International Journal}, 13(2-3),
  1995.

\bibitem[Mal95]{malgrangeexpositiones}
B.~Malgrange.
\newblock Sommation des s\'eries divergentes.
\newblock {\em Expositiones Mathematicae. International Journal},
  13(2-3):163--222, 1995.

\bibitem[MZ00]{MarotteZhang}
F.~Marotte and C.~Zhang.
\newblock Multisommabilit\'e des s\'eries enti\`eres solutions formelles d'une
  \'equation aux {$q$}-diff\'erences lin\'eaire analytique.
\newblock {\em Annales de l'Institut Fourier}, 50(6):1859--1890, 2000.

\bibitem[Ram92]{ramistoulouse}
J.-P.~Ramis.
\newblock About the growth of entire functions solutions of linear algebraic
  {$q$}-difference equations.
\newblock {\em Toulouse. Facult\'e des Sciences. Annales. Math\'ematiques.
  S\'erie 6}, 1(1):53--94, 1992.

\bibitem[Ram93]{RamisDivergentes}
J.-P.~Ramis.
\newblock {\em S\'eries divergentes et th\'eories asymptotiques}, volume 121 of
  {\em Panoramas et Synth\`ese, suppl.}
\newblock SMF, 1993.

\bibitem[RM90]{raminet}
J.-P.~Ramis and J.~Martinet.
\newblock Th\'eorie de {G}alois diff\'erentielle et resommation.
\newblock In {\em Computer algebra and differential equations}, Comput. Math.
  Appl., pages 117--214. Academic Press, 1990.

\bibitem[RZ02]{RZ}
J.-P.~Ramis and C.~Zhang.
\newblock D\'eveloppement asymptotique $q$-{G}evrey et fonction th\^eta de
  {J}acobi.
\newblock {\em C. R., Math., Acad. Sci. Paris}, 335(11):899--902, 2002.

\bibitem[Sau00]{Sfourier}
J.~Sauloy.
\newblock Syst\`emes aux $q$-diff\'erences singuliers r\'eguliers:
  classification, matrice de connexion et monodromie.
\newblock {\em Annales de l'Institut Fourier}, 50(4):1021--1071, 2000.

\bibitem[Sau02]{Spoligono}
J.~Sauloy.
\newblock La filtration canonique par les pentes d'un module aux
  $q$-diff\'{e}rences.
\newblock Pr\'epublication du Laboratoire Emile Picard n.249. {\tt
  arXiv:math.QA/0210221}, 2002.

\bibitem[Sau04]{sauloyfiltration}
J.~Sauloy.
\newblock La filtration canonique par les pentes d'un module aux
  {$q$}-diff\'erences et le gradu\'e associ\'e.
\newblock {\em Annales de l'Institut Fourier}, 54(1):181--210, 2004.

\bibitem[WW88]{WW}
E.T.~Whittaker and G.N.~Watson.
\newblock {\em A Course of Modern Analysis}.
\newblock Cambridge University Press, Cambridge, 1988.

\bibitem[Zha99]{zhangfourier}
C.~Zhang.
\newblock D\'eveloppements asymptotiques {$q$}-{G}evrey et s\'eries
  {$Gq$}-sommables.
\newblock {\em Annales de l'Institut Fourier}, 49(1):227--261, 1999.

\bibitem[Zha00]{zhangthetafunction}
C.~Zhang.
\newblock Transformations de {$q$}-{B}orel-{L}aplace au moyen de la fonction
  th\^eta de {J}acobi.
\newblock {\em Comptes Rendus de l'Acad\'emie des Sciences. S\'erie I.
  Math\'ematique}, 331(1):31--34, 2000.

\bibitem[Zha01]{zhangqGamma}
C.~Zhang.
\newblock Sur la fonction {$q$}-gamma de {J}ackson.
\newblock {\em Aequationes Mathematicae}, 62(1-2):60--78, 2001.

\bibitem[Zha02]{zhanggroningen}
C.~Zhang.
\newblock Une sommation discr\`ete pour des \'equations aux {$q$}-diff\'erences
  lin\'eaires et \`a coefficients analytiques: th\'eorie g\'en\'erale et
  exemples.
\newblock In {\em Differential equations and the Stokes phenomenon}, pages
  309--329. World Sci. Publishing, River Edge, NJ, 2002.

\bibitem[Zha03]{zh2}
C.~Zhang.
\newblock Sur les fonctions {$q$}-{B}essel de {J}ackson.
\newblock {\em Journal of Approximation Theory}, 122(2):208--223, 2003.

\end{thebibliography}
%%%%%%%%%%%%%%%%%%%%%%%%%%%%%%%%%%%%%%%%%%%%%%%%%%%%%%%%%%%%%%%%%%%%%%%%%%%%%%%%%%%%%%%%
\newcommand{\noopsort}[1]{}

\end{document}